\documentclass[11pt]{article}
\usepackage{amsfonts}
 \usepackage{color}
\usepackage{graphics}

\usepackage{cite}
\usepackage{latexsym}
\usepackage{amsmath}
\usepackage{amssymb}
\usepackage[dvips]{epsfig}
\usepackage{amscd}

\usepackage[paper=a4paper, left=2.1cm, right=2.1cm, top=2.2cm, bottom=1.5cm, headheight=5.5pt, footskip=0.8cm, footnotesep=0.8cm, centering, includefoot]{geometry}

\newtheorem{theorem}{Theorem}[section]
\newtheorem{remark}{Remark}[section]

\newtheorem{definition}{Definition}[section]
\newtheorem{lemma}[theorem]{Lemma}

\newcommand{\rt}{\rightarrow}

\newcommand{\n}{\rho}
  \newcommand{\nv}{\rho_\ve}
   \newcommand{\uv}{u_\ve}
  
  \newcommand{\ioo}{\int_0^T\int }
  
\newcommand{\ti}{\tilde}

\newcommand{\ro}{\rightarrow}

\renewcommand{\div}{ {\rm div }  }

\newcommand{\vf}{v_\ve}

\renewcommand{\r}{\mathbb{R}}

\renewcommand{\b}{Q_{\ve }}

\newcommand{\bt}{\begin{theorem}}
\newcommand{\bl}{\begin{lemma}}
\newcommand{\el}{\end{lemma}}
\newcommand{\et}{\end{theorem}}
\newcommand{\ga}{\gamma}

\newcommand{\ve}{\varepsilon}
\newcommand{\la}{\label}

 \newcommand{\ka}{\kappa}

\newcommand{\bn}{\begin{eqnarray}}
\newcommand{\en}{\end{eqnarray}}
\newcommand{\bnn}{\begin{eqnarray*}}
\newcommand{\enn}{\end{eqnarray*}}

\newcommand{\bnnn}{\begin{eqnarray*}}
\newcommand{\ennn}{\end{eqnarray*}}

\newcommand{\ba}{\begin{aligned}}
\newcommand{\ea}{\end{aligned}}
\newcommand{\be}{\begin{equation}}
\newcommand{\ee}{\end{equation}}
\def\O{{\Omega }}

\def\p{\partial}
\def\norm[#1]#2{\|#2\|_{#1}}

\newcommand{\si}{\sigma}

\def\la{\label}

\def\na{\nabla}

\def\xl{\left}
\def\xr{\right}

\def\xd{\cdot}

\makeatletter      
\@addtoreset{equation}{section}
\makeatother

\makeatletter      
\@addtoreset{equation}{section}
\makeatother       

\title{Global existence of weak solutions to the compressible quantum Navier-Stokes equations with degenerate viscosity
}
 \author{ Boqiang L\"u\thanks{College of Mathematics and Information Science, Nanchang Hangkong University, Nanchang 330063, People's Republic of China ({\tt lvbq86@163.com}). }
 \quad Rong Zhang\thanks{The Institute of Mathematical Sciences, The Chinese University of Hong Kong, Shatin, Hong Kong ({\tt rzhang@math.cuhk.edu.hk}).}
 \quad Xin Zhong\thanks{School of Mathematics and Statistics, Southwest University, Chongqing 400715, People's Republic of China ({\tt xzhong1014@amss.ac.cn}). }}

\date{}

\begin{document}
 \maketitle

\begin{abstract}
We study the compressible quantum Navier-Stokes (QNS) equations with degenerate viscosity  in the three dimensional periodic domains. On the one hand, we consider  QNS with     additional damping terms. Motivated by the recent works [Li-Xin, arXiv:1504.06826] and [Antonelli-Spirito, Arch. Ration. Mech. Anal., 203(2012), 499--527], we construct a suitable approximate system which has smooth solutions satisfying the energy inequality and the BD entropy estimate. Using this system,
we obtain the global existence of weak solutions to the compressible QNS equations  with damping terms for large initial data. Moreover,   we obtain some new a priori estimates,  which can avoid using the assumption that the gradient of the velocity is a well-defined function,  which is indeed used directly in  [Vasseur-Yu, SIAM J. Math. Anal., 48 (2016), 1489--1511; Invent. Math.,  206 (2016), 935--974].
On the other hand, in the absence of damping terms, we also prove the  global existence of weak solutions to the compressible QNS equations without the  lower bound assumption  on the dispersive coefficient, which improves the previous result due to   [Antonelli-Spirito, Arch. Ration. Mech. Anal., 203(2012), 499--527].
\end{abstract}

Keywords: compressible quantum Navier-Stokes equations; global weak solutions; degenerate viscosities; vacuum.

Math Subject Classification: 35Q35; 76N10

\section{Introduction}
The  quantum Navier-Stokes equations with damping terms which read as follows:
\be\la{qns}
\begin{cases}
\n_t+\div(\n u)  = 0, \\
(\n u)_t+ \div(\n u \otimes u)-2\nu\div(\n\mathcal{D} u) +\na P -2\ka^2\n\na\left(\frac{\Delta \sqrt{\n}}{\sqrt{\n}}\right)+r_0u+r_1\n|u|^2u=0.
\end{cases}
\ee
Here, $x\in \Omega \subset \mathbb{R}^3$,  $t>0,$ $\n$ is the density, $u=(u_1,u_2,u_3)$ is the velocity field, $ \mathcal{D}u =\frac{1}{2}(\na u+(\na u)^{\rm tr} )$ is the symmetric part of the velocity gradient,
$P(\rho)=a\rho^\ga (a>0,\ga> 1)$ is the pressure. Without   loss of generality, it is  assumed that $a=1$. The positive constants $\nu$ and $\kappa$ are the viscosity and the dispersive coefficients, respectively. The constants $r_0$ and $r_1$ in the damping terms are all positive.
Let $\O=\mathbb{T}^3$   be the three dimensional torus, we consider the system \eqref{qns}  with periodic boundary conditions. The initial conditions are imposed as
\begin{equation}\label{1.2}
\rho(x,0)=\rho_0(x),\ \ \rho u(x,0)=m_0(x).
\end{equation}

When $r_0=r_1=0$, i.e., there is no damping terms, the   system \eqref{qns} is a special case of the Navier-Stokes-Korteweg (NSK)  equations, which reads as
\be\la{qns1}
\begin{cases}
\n_t+\div(\n u)  = 0, \\
(\n u)_t+ \div(\n u \otimes u) +\na P- \div\mathbb{S} -\div\mathbb{K}=0.
\end{cases}
\ee
The viscosity stress tensor $\mathbb{S}$ and the capillarity (dispersive) term $\mathbb{K}$ are defined by
\be\la{aa'1} \ba \mathbb{S}\triangleq h \mathcal{D} u+g\div u\mathbb{I} \ea\ee
and \be\la{naa'1} \ba \mathbb{K}\triangleq  \xl(\n\div(k(\n)\na \n-\frac{1}{2}(\n k'(\n)-k(\n))|\na \n|^2\xr)\mathbb{I}-k(\n)\na \n\otimes\na\n ,\ea\ee
where $\mathbb{I}$ is the identical matrix, and $h,g$ satisfy the physical
restrictions \be\notag\la{bd2} h>0,\quad h+3g\ge 0.\ee
Indeed, choosing
\be
\la{hgvv1}h(\n)=2\nu\n,\quad g(\n)=0,\quad k(\n)=\frac{\kappa^2}{\n},
\ee
the NSK equations \eqref{qns1} becomes the QNS one \eqref{qns} without damping terms. For more detailed  derivation of the QNS equations, please refer to \cite{JU12}.
 In particular, the QNS equations without viscosity ($\nu=0$) is the    Quantum Hydrodynamics model for superfluids (see \cite{LL77}), whose global weak solutions with finite energy was studied in \cite{AM09,AM12}.
It is well known that the  NSK equations reduces to the Navier-Stokes (NS) equations when there is no capillarity (dispersive) term $\mathbb{K}$.
 One of the main difficulties in studying the compressible NS (or QNS, NSK) equations
 with degenerate viscosity coefficients is to
   estimate  the gradient of the velocity field in the vacuum region,
please refer to
\cite{AM09,AM12,AS2017,AS2015,BreschDesjardins02,
BreschDesjardins03,
BreschDesjardins03b,BreschDesjardins05,
BreschDesjardinsLin05,BD06jmpa,fe1,GL2015,gjx1,glx,hl2,hl1,JiangXinZhang06,jwx1,jun,jun2011,LV18jmpa,llx06,li04,Lions98,LiuXinYang98,lvp,MelletVasseur05,vy1,vy,Yangzhu02} and the references therein.

For the one dimensional space,  the global existence of weak solutions for the QNS equations was proved  by J\"ungel  \cite{jun2011}. Then, for weak solutions   required a special choice of the test function   $\n \phi$ with $\phi$ smooth and compactly supported, he \cite{jun} also obtained the global weak solutions to the three dimensional  QNS equations in the case $\kappa>\nu$ and $\gamma>3.$
Very recently, for $\ga,\ka,$ and $\nu$ satisfying  \be\la{asrq} \begin{cases} 1<\gamma,~~\kappa<\nu,~~~~&\O=\mathbb{T}^2,\\1<\gamma<3,~~ \kappa<\nu<\frac{3\sqrt{2}}{4}\kappa,~~~~&\O=\mathbb{T}^3, \end{cases}  \ee Antonelli-Spirito \cite{AS2017}  proved the global existence of finite energy weak solutions, which  is   the first result of global existence for finite energy weak solutions to NSK equations in high dimensional space.
As mentioned in \cite{AS2017}, one of the key ideas in \cite{AS2017} is to construct     proper smooth  approximating solutions, which is motivated by the parabolic regularization methods owing to  Li-Xin \cite{li04}. Indeed,   Li-Xin \cite{li04}     proved  the global existence of finite energy weak solutions to the compressible NS  equations with general degenerate viscosity coefficients in two or three dimensional periodic domains or whole spaces, which in   particular   solved an open problem proposed by Lions \cite{Lions98}.

%
Furthermore,
there are many works considering  the compressible NS (or QNS, NSK) equations by considering the system with some additional terms, such as a cold pressure term, the damping terms or other source terms (please see \cite{BreschDesjardins03,BD06jmpa,BreschDesjardinsLin05,GL2015,vy1,LV18jmpa}  and the references therein).
In particular, Vasseur-Yu \cite{vy1}  considered  global existence of finite energy weak solutions of the QNS equations with damping terms \eqref{qns}. 
 Then,
using  the global weak solutions to system \eqref{qns} obtained in  \cite{vy1}, by different methods from those in Li-Xin \cite{li04}, Vasseur-Yu \cite{vy}  studied  the  global weak solutions to the compressible NS equations \eqref{qns1}-\eqref{hgvv1} with $\ka=0$.
%
 The key issues in Vasseur-Yu \cite{vy1,vy} rely crucially  on the assumptions that   $\na u$ is a well-defined function and that $\sqrt{\n}\na u \in L^2((0,T)\times\O)$,   which are confused for us (see Remark \ref{key1} below for more details).
Indeed, it seems impossible to define $\na u$ as    functions
without enough  regularity of $u$ due to the high degenerate viscosity at   vacuum.
Hence,  in this paper, we will  reconstruct  suitable   approximate system   to obtain  the  global existence of weak solutions to system  \eqref{qns}. Moreover, the weak solutions are more regular than   those obtained by  Vasseur-Yu \cite{vy1}
and can be used to obtain the global weak solutions to the compressible NS equations with degenerate viscosity. This will be shown in a forthcoming paper \cite{lvp}.
Furthermore, we also  improve  the restriction on the range of $\kappa$ in \cite{AS2017} by removing the lower bound   $\frac{4}{3\sqrt{2}}\nu $.

Now, we explain the notations and conventions used throughout this paper.  For  $\O=\mathbb{T}^3$, set
$$  \int \cdot ~dx\triangleq\int_{\O}\cdot ~dx.$$ Moreover, for $1\le r\le \infty$ and  $k\ge 1, $  the standard Lebesgue and  Sobolev spaces are defined as follows:
   \bnn
L^r=L^r(\O ),\quad
W^{k,r}  = W^{k,r}(\O) , \quad H^k = W^{k,2}.
 \enn
We will consider the  problem \eqref{qns}--\eqref{1.2} with  the initial data $\n_0,m_0$ satisfying that
\be\la{pini1}
\begin{cases} \n_0\ge 0 \mbox{ a.e. in }\O,~~~\n_0\not\equiv 0,~~~\n_0 \in L^1\cap L^\ga,~~~ \na\sqrt{\n_0}\in L^2,\\
 m_0\in L^{1}, \, m_0= 0 \mbox{ a.e. on }\O_0,~~~
 \n_0^{-1} m_0^2\in L^{1},\\
-r_0\log_{-}\n_0 \in L^1,  ~\mbox{with}~\log_{-}g\triangleq \log \min\{1,g\},
 \end{cases}
 \ee
 where $\O_0$ is the vacuum set of $\n_0 ,$   defined by \be \la{oox1}\O_0\triangleq\{x\in \O\,|\n_0(x)=0\}.\ee

 Next,
we give the definition of a  weak solution to \eqref{qns}--\eqref{1.2}.
\begin{definition}\la{def1} Let $\O=\mathbb{T}^3$,  $(\n,u)$ is said to be a weak solution to \eqref{qns}--\eqref{1.2} if
\be \ba\label{ims1}\begin{cases}0\le \n\in L^\infty(0,T;L^1\cap L^\gamma ),
\\ \na\n^{\frac{\ga}{2}}\in L^2(0,T;L^2 ),
\\ \na\sqrt{\n},\,\,\sqrt{\n}u\in L^\infty(0,T;L^2 ),
\\ \na( \sqrt{\n}u),~~\na( \sqrt{\n}u)-u\otimes \na \sqrt{\n}~\in L^2(0,T;L^2 ),
\\\sqrt{r_0}u ~\in L^2(0,T;L^2 ),~~r_1^{1/4}\n^{1/4}u ~\in L^4(0,T;L^4 ),
 \\  \kappa \na^2 \sqrt{\n} ,~\kappa\sqrt{\n}\na^2\log \n \in L^2(0,T;L^2 ),
\end{cases}\ea\ee
with $ (\n,\sqrt{\n}u)$ satisfying \be\la{fin1} \begin{cases} \n_t+\div (\sqrt{\n}\sqrt{\n}u)=0,\\ \n(x,t=0)=\n_0(x),\end{cases} \mbox{ in }\mathcal{D}', \ee
and if the following equality holds for all smooth test function $\phi(x,t)$ with compact support such that $\phi(x,T)=0:$
\be \ba\label{fin2}
&\int  m_0\cdot \phi(x,0)dx+\int_0^T \int  \left(\sqrt{\n}\sqrt{\n}u\cdot\phi_t+\sqrt{\n}u\otimes \sqrt{\n}u:\na\phi + \n^\ga \div \phi \right)dxdt
\\& \quad - \nu\int_0^T \int \sqrt{\n}\xl( (\na(\sqrt{\n}u)-u\otimes \na\sqrt{\n}) :\na \phi +(\na^{tr}(\sqrt{\n}u)-\na\sqrt{\n}\otimes u) :\na \phi \xr)dxdt \\
& =\int_0^T \int \xl(r_0 u\cdot\phi  +r_1  \n|u|^2u\cdot\phi +4\kappa^2  \Delta \sqrt{\n}\na \sqrt{\n} \cdot  \phi +2\kappa^2 \Delta \sqrt{\n}  \sqrt{\n}\div\phi\xr) dxdt.
\ea\ee
\end{definition}

Our first result  reads as follows:
\begin{theorem}\la{qvth2}
Suppose that $\ga\in (1,3)$
and  $11\kappa \le  \nu  $. 
Moreover, assume that the initial data $(\n_0,m_0)$ satisfy \eqref{pini1}. Then, there exists a global weak solution $(\n,u)$ to the problem \eqref{qns}--\eqref{1.2} satisfying
\be\label{19ims}\ba& \sup_{0\le t\le T}\int\left(\n |u|^2 +  \n^\ga \right)dx  + \int_0^T\int  \xl( r_0 |u|^2  +r_1 \n |u|^4\xr) dxdt  \le C,\ea\ee
\be\label{19ims2}\ba& \sup_{0\le t\le T}\int |\na \sqrt{\n}|^2 dx    +\int_0^T \int \left(|\na(\sqrt{\n}u)-u\otimes \na\sqrt{\n}|^2 + |\na \n^{\frac{\gamma}{2}}|^2\right)dxdt\le C+Cr_0+Cr_1,  \ea\ee
and
\be\label{19ims3}\ba \kappa^2\int_0^T\int\xl(|\nabla\rho^\frac14|^4+|\nabla^2\rho^\frac12|^2\xr)dxdt
+r_1\kappa\int_0^T \int  |\na (\sqrt{\n} u )|^2  dxdt \le C+Cr_0+Cr_1.  \ea\ee
where  $C$ is a positive generic constant depending only on the initial data, but
independent of $\kappa, r_0$, and $r_1$.
\end{theorem}

 A few remarks are in order:

\begin{remark}\label{key1}
It should be noted that  the arguments in Vasseur-Yu \cite{vy,vy1}
rely crucially  on  the assumption that the gradient of  velocity  field  $\na u$ is a well-defined  function, which indeed does not make sense in the presence of vacuum.
%
 In particular, in the proof of  \cite[Lemma 4.2]{vy}, which is crucial to deduce the key Mellet-Vasseur type estimate  in \cite{vy}, it requires  essentially that $ \na u$  is a  well-defined function.

Very recently, Lacroix-Violet \& Vasseur \cite{LV18jmpa}  also study  the QNS equations and   consider  a    new  function  $\mathbb{T_\nu}\in L^2(\O\times (0,T))$ satisfying 
\be \label{nims4}  \sqrt{\nu\n}\mathbb{T_\nu}=\nu\na (\n u)-2\nu\sqrt{\n}u\otimes\na \sqrt{\n}. \ee
More precisely, they \cite{LV18jmpa}  use  the  function $\mathbb{T_\nu}$ to give a new understanding of $\sqrt{\n} \na u$. However,  as mentioned in  \cite{LV18jmpa},  it still does not allow  to define  the  gradient of  velocity $ \na u$  as a function. \end{remark}

\begin{remark}\label{key2} If $\ka>0$ and $r_1>0,$ Theorem \ref{qvth2} shows that $\sqrt{\n}u\in L^2(0,T;H^1(\O))$, which is a complete  new regularity estimate.  Combining this fact with  $ \sqrt{\n} \in L^2(0,T; H^1(\O))$ shows that \be \na(\n u)=\sqrt{\n}\na(\sqrt{\n} u)+\na\sqrt{\n}\otimes \sqrt{\n}u,\ee holds rigorously in the sense of function.   This new observation is   helpful for further studies on   the weak solutions of compressible Navier-Stokes equations, which will be shown in our another  paper \cite{lvp}.
\end{remark}



Next, we also obtain the global weak solutions to  system \eqref{qns} without damping terms.
\begin{theorem}\la{newth}
 Suppose that $r_0=r_1=0$, $\ga\in (1,3)$,
and  $11\kappa \le  \nu   $. 
Moreover, assume that the initial data $(\n_0,m_0)$ satisfy \eqref{pini1}$_1$,  \eqref{pini1}$_2$, and
 \begin{align}\label{MVinitial}
 \sqrt{\rho_0} \in L^{2+\eta},~~~ \sqrt{\rho_0 }u_0 \in L^{2+\eta},
 \end{align}
 for any $\eta>0$. Then the problem \eqref{qns}--\eqref{1.2} admits  a global weak solution $(\n,u)$ satisfying \eqref{ims1}$_1$--\eqref{ims1}$_3$. Moreover,   $(\n,\sqrt{\n}u)$ satisfy  \eqref{fin1} and
 \be \ba\label{nfin2}
&\int  m_0\cdot \phi(x,0)dx+\int_0^T \int \left(\sqrt{\n}\sqrt{\n}u\cdot\phi_t+\sqrt{\n}u\otimes \sqrt{\n}u:\na\phi + \n^\ga \div \phi \right)dxdt
\\&  - 2\nu\int_0^T \int (\sqrt{\n}u\otimes\na\sqrt{\n}) :\na \phi dxdt - 2\nu\int \int_0^T (\na\sqrt{\n}\otimes\sqrt{\n}u) :\na \phi dxdt \\
 &  +\nu\int_0^T \int \sqrt{\n}\sqrt{\n}u\cdot\Delta \phi dxdt +\nu\int_0^T \int  \sqrt{\n}\sqrt{\n}u\cdot\na\div \phi  dxdt \\
&-4\kappa^2\int_0^T \int(\na \sqrt{\n}\otimes\na\sqrt{\n}): \na \phi dxdt +2\kappa^2\int_0^T \int\sqrt{\n}\na\sqrt{\n}\cdot\na\div\phi dxdt=0,
\ea\ee
where $\phi(x,t)$ is  a smooth test function with compact support satisfying $\phi(x,T)=0$.
 \end{theorem}

\begin{remark} Compared with \cite{AS2017},  our Theorem 1.2 succeeds in removing their assumption on  the lower bound of dispersive coefficient  $ \kappa>\frac{4}{3\sqrt{2}}\nu.$
\end{remark}

We now sketch some main ideas used in our analysis.
The main point of this paper is to construct smooth approximate solutions
satisfying the energy inequality and the BD entropy estimate.
%
%
Thanks to Li-Xin \cite{li04}, we first propose to approximate \eqref{qns}$_1$ by a parabolic equation \eqref{1.12}$_1$.
Next, on the one hand, some similar regularization  in \eqref{1.12}$_2$ as those in \cite{li04} are  considered accordingly with respect to the parabolic  regularization in \eqref{1.12}$_1$. On the other hand, the third order  capillarity  term will bring us some new difficulties. Motivated by \cite{AS2017,AS2015}(see also \cite{jun}), by using  the effective velocity $w\triangleq u+\mu\na\log\n$ with $\mu\triangleq\nu-\sqrt{\nu^2-\kappa^2}$ to handle the third order  capillarity  term, we thus need some additional regularization terms of $\na\log\n$ in \eqref{1.12}$_2$.  As a result,  we consider the following  approximate system
 \be \label{1.12}
 \begin{cases}
 \n_t+\div(\n u)    = \ve v\div( |\na v|^2 \na v)+\ve\n^{-p_0},
 \\   \n u_t+ \n u\cdot \na u   -2\nu\div( \n    \mathcal{D} u) +\na P -2\ka^2\n\na\left(\frac{\Delta v}{v}\right)+r_0u+r_1\n|u|^2u\\=  \sqrt{\ve}\div(\n \na u)+\sqrt{\ve}\mu\div(\n \na^2\log \n) +\ve v |\na v|^2 \na v\cdot \na u+\ve\mu v |\na v|^2 \na v\cdot \na (\na \log\n)\\ \quad -\ve \n^{-p_0} u- \ve^{\frac32}\n    |w|^{3}u-\ve\mu\na \n^{-p_0} - \ve\mu\na( v\div( |\na v|^2 \na v))+\ve\mu v\div( |\na v|^2 \na v)\na \log \n,
 \end{cases}
 \ee
where $v\triangleq\sqrt{\n}$.
First of all, following the similar arguments as those in \cite{li04},   the smooth solutions to the approximate system \eqref{1.12} satisfy both the energy inequality and the BD entropy estimates (see \eqref{zqbb30} and \eqref{zqbb5'}). Then, we use a De Giorgi-type procedure to bound the density from above and below (see \eqref{qc1}), provided the initial density is strictly away from vacuum. In particular, it is proved   that  the density is strictly away from vacuum.

With these estimates in hand, we will dedcue the higher order estimates on $(\n,u)$,  which are necessary to get the global strong solutions to the system \eqref{1.12}. However, due to the third order  capillarity  term, it is difficult to establish  directly the desired higher order estimates on $(\n,u)$. To this end, we consider the solutions $(\n,w)$ to a transformation system \eqref{w=u0}, which is equivalent to the system \eqref{1.12} of $(\n,u)$.  Then,   by using the $L^p$-theory for parabolic equations, we  get the desired estimates on  $(\n,w)$ and thus the estimates on  $(\n,u)$ (see \eqref{hi0} and \eqref{zqbc4}).
 This  implies that  the approximate system \eqref{1.12} has a global strong solution with smooth initial data.
Next, after adapting the compactness results due to  \cite{BreschDesjardins03,BreschDesjardinsLin05,AS2015},
we can  obtain the global existence of the weak solutions to \eqref{qns}  and thus prove Theorem 1.1.

Finally, for the system \eqref{qns} without damping terms, we will consider the approximate system \eqref{1.12} with $ r_0=r_1=0 $.
In the absence of damping terms, we need further to derive the Mellet-Vasseur type estimate. As pointed in \cite{AS2015,AS2017,vy}, 
the third order dispersive term prevents one from obtaining directly a Mellet-Vasseur type inequality. This difficulty is overcome by deriving the Mellet-Vasseur type estimate on  $(\n,w)$  to the transformation system \eqref{w=u0} without third order term. Therefore, it  shows that the approximate system \eqref{1.12} with $ r_0=r_1=0$   has smooth solutions satisfying the
energy inequality, the BD entropy one, and the Mellet-Vasseur type estimate. The compactness results  \cite{AS2015} ensure Theorem 1.2 directly.

The rest of the paper is organized as follows. In Section 2, we construct the approximate system and derive the a priori estimates. Section 3 is devoted to compactness results of the approximate solutions. Theorem \ref{qvth2} is proved in Section 4. Finally, Section 5 will show the Mellet-Vasseur type inequality to the system \eqref{qns} without damping terms and  then prove Theorem \ref{newth}.


\section{A priori estimates}

Let $v\triangleq\n^{1/2}$ and
\be \label{lw} w\triangleq u+\mu\na \log \n \ee
with $\mu=\nu-\sqrt{\nu^2-\kappa^2}$ and  $11\kappa \le  \nu  $, we  consider  the following approximate system
 \be\label{zqba1} \begin{cases}   \n_t+\div(\n u)    = \ve v\div( |\na v|^2 \na v)+\ve\n^{-p_0},
     \\   \n u_t+ \n u\cdot \na u   -2\nu\div( \n    \mathcal{D} u) +\na P -2\ka^2\n\na\left(\frac{\Delta v}{v}\right)+r_0u+r_1\n|u|^2u\\=  \sqrt{\ve}\div(\n \na u)+\sqrt{\ve}\mu\div(\n \na^2\log \n) +\ve v |\na v|^2 \na v\cdot \na u+\ve\mu v |\na v|^2 \na v\cdot \na (\na \log\n)\\ \quad -\ve \n^{-p_0} u- \ve^{3/2}\n    |w|^{3}u -\ve\mu\na \n^{-p_0}- \ve\mu\na( v\div( |\na v|^2 \na v))+\ve\mu v\div( |\na v|^2 \na v)\na \log \n,  \end{cases}\ee
  where the constants $p_0 $ and $\ve$
  satisfying
  \be\notag\la{qq1}p_0=50,\quad 0<\ve\le  10^{-10}.\ee
%
  The initial conditions of the system \eqref{zqba1} are imposed as:
\be \label{zqbb35} (\n,u)(x,0)= (\n_{0\ve},u_{0\ve}),\ee
where  smooth  $\O$-periodic functions $ \n_{0\ve}>0$ and $u_{0\ve} $ satisfying
\be \la{qpd9}\ba&\|r_0\log_{-}\n_{0\ve}\|_{L^1}+\|\n_{0\ve}\|_{L^1\cap L^{\ga}} +\|\na  \n_{0\ve}^{ 1/2 } \|_{L^2}+\ve \|\na  \n_{0\ve}^{1/2 } \|^4_{L^4}+\ve \|  \n_{0\ve}^{-p_0 } \|_{L^1} \le C \ea \ee
and
\be  \la{qpd22}\int  \n_{0\ve}   |u_{0\ve}|^{2}dx\le C\ee
for some constant $C$ independent of $\ve.$

 Some alternative ways of the third order tensor term are stated as follows
  \be \ba\label{len0}  2 \n\na\left(\frac{\Delta v}{v}\right)=\div(\n \na^2\log \n) =\na \Delta\n-4\div(\na v\otimes\na v).\ea\ee
 Let $T>0$ be a fixed time and $(\rho, u)$  be
a smooth solution to \eqref{zqba1}--\eqref{zqbb35} on
$\O \times (0,T]. $
Then, we will establish some necessary a priori bounds
for  $(\rho,u)$. The first one is the energy-type inequality.
\begin{lemma} \label{zqlem10} Suppose that   $11\kappa \le  \nu  $, then there exists some generic constant $C$ independent of $\ve$, $r_0$, $r_1$, and $\kappa$  such that
\be\label{zqbb30}\ba& \sup_{0\le t\le T}\int\left(\n |u|^2+\n+  \n^\ga+\ve\n^{-p_0} +(2\kappa^2+2\mu\sqrt{\ve})|\na v|^2+\ve\mu|\na v|^4 \right)dx  \\
 &+ \nu \int_0^T\int    \n   |\mathcal{D} u|^2 dxdt + r_0\int_0^T\int    |u|^2 dxdt+r_1\int_0^T\int    \n   |u|^4 dxdt\\&+ \sqrt{\ve}\int_0^T\int    \n   |\na u|^2 dxdt +\ve\int_0^T\int   \left(   |\na v|^4 + |\na v|^4|u|^2   + \n^{-p_0}|u|^2+\ve^{1/2}\n  |w|^{3}|u|^{2}\right) dxdt\\
 &+(2\kappa^2+2\mu\sqrt{\ve})\ve\int_0^T\int   \left(   |\na v|^2|\na^2 v|^2  +|\na|\na v|^2 |^2 + (2p_0+1)|\na v|^2v^{-2p_0-1}\right) dxdt\\ & +\ve^2 \int_0^T    \int \left(\mu|\na v|^4|\na^2 v|^2+\mu|\na v|^4|\na|\na v|  |^2
 +(2p_0+1)\nu|\na v|^4v^{-2p_0-2}+\n^{-2p_0-1}\right)dx  dt  \le C.\ea\ee
\end{lemma}
{\it Proof.} First, integrating \eqref{zqba1}$_1$ over $\O$ yields
\be \label{zqbc6} \xl(\int\n dx\xr)_t +\ve \int  |\nabla v|^4dx  =\ve \int \n^{-p_0}dx .\ee

Next, multiplying $\eqref{zqba1}_2$ by $u $ and  integrating the resulting equations by parts, we obtain after using $\eqref{zqba1}_1 $ that
\be\label{zqbb61}\ba& \frac12 \xl(\int\n |u|^2dx\xr)_t   +2\nu\int \n  |\mathcal{D}u|^2  dx +\sqrt{\ve}\int \n  |\na u|^2  dx
+\frac{\ve}{2}\int \n^{-p_0} |u|^2dx\\
&\quad+\ve^{3/2}\int\n  |w|^{3}|u|^{2}dx +r_0\int |u|^{2}dx+r_1\int\n  |u|^{4}dx+\int   u \cdot \na\n^\ga  dx+\ve\mu\int   u \cdot \na\n^{-p_0}  dx\\
&= \frac{\ve}{2}\int  v\div (|\na v|^2\na v) |u|^2dx +\ve\int v |\na v|^2 \na v\cdot \na u\cdot udx \\
&\quad+\ve\mu \int v |\na v|^2 \na v\cdot \na (\na \log\n)\cdot udx+2(\kappa^2+\sqrt{\ve}\mu)\int \n \na\left(\frac{\Delta v}{v}\right)\cdot udx \\
&\quad +\ve\mu \int v\div( |\na v|^2 \na v)\na \log \n \cdot udx- \ve\mu\int\na( v\div( |\na v|^2 \na v))\cdot udx\\
&=\sum_{i=1}^6I_i.\ea\ee

Integration by parts gives
\be\ba\label{len1}
I_1+I_2
&=   -\frac{\ve}{2}\int   |\na v|^4 |u|^2dx.\ea\ee

Since $\na \log\n=2v^{-1}\na v$, one has
\be\ba\label{len1'}
 I_3&= \ve \mu\int   v |\na v|^2 \na v\cdot \na (\na \log\n)\cdot udx \\
 &= 2\ve \mu\int   |\na v|^2 \na v\cdot \na^2v \cdot udx-2\ve \mu\int   v^{-1}|\na v|^4 \na v\cdot udx \\
 &\le \frac{\ve}{4}\int   |\na v|^4|u|^2dx
 +8\ve\mu^2\int   |\na v|^2|\na^2v|^2dx+8\ve\mu^2\int v^{-2}|\na v|^6   dx\\
 &\le  \frac{\ve}{4}\int   |\na v|^4|u|^2dx+24\ve\mu^2\int   |\na v|^2|\na^2v|^2dx+64\ve\mu^2\int|\na |\na v|^2 |^2dx,
 \ea\ee
where in the last inequality one has used the following fact
\be\ba \label{len1-5}
 \int v^{-2}|\na v|^6   dx \le  2\int|\na v|^2 |\Delta v|^2dx +8\int|\na |\na v|^2 |^2dx.
 \ea\ee
Indeed,  integration by parts together with some directly calculations show that
 \be\ba \label{len1-3}
 \int v^{-2}|\na v|^6   dx&=  \int v^{-2}|\na v|^4   \na v\cdot\na vdx \\
  &=-\int v\na v^{-2} |\na v|^4   \cdot\na vdx-\int v^{-1}\na\ |\na v|^4 \cdot\na vdx- \int v^{-1}|\na v|^4 \Delta vdx\\
    &= 2\int  v^{-2} |\na v|^6dx- 2\int v^{-1}|\na v|^2\na|\na v|^2 \cdot\na vdx- \int v^{-1}|\na v|^4 \Delta vdx,
 \ea\ee
 that is
 \be\ba \label{len1-4}
 \int v^{-2} |\na v|^6 dx&=2\int v^{-1}|\na v|^2\na|\na v|^2 \cdot\na vdx+ \int v^{-1}|\na v|^4 \Delta vdx\\
 &\le \frac{1}{2}\int v^{-2} |\na v|^6 dx+  \int|\na v|^2 |\Delta v|^2dx +4\int|\na |\na v|^2 |^2dx.
 \ea\ee
This yields \eqref{len1-5} directly.

For the term $I_4$, it deduces from \eqref{zqba1}$_1$ and integration by parts that
\be\ba\label{len4}
 I_4&
 =-2(\kappa^2+\sqrt{\ve}\mu)\int  \frac{\Delta v}{v}\div(\n u)dx\\
 &= -2(\kappa^2+\sqrt{\ve}\mu)\int  \frac{\Delta v}{v}\left(-2vv_t +  \ve v\div( |\na v|^2 \na v)+\ve\n^{-p_0}\right) dx\\
&= -2(\kappa^2+\sqrt{\ve}\mu) \frac{d}{dt}\int |\na v|^2dx\\
&\quad-2(\kappa^2+\sqrt{\ve}\mu)\ve\int \xl(|\na v|^2 |\na^2 v|^2+\frac{1}{2}|\na |\na v|^2|^2+(2p_0+1)|\na v|^2v^{-2p_0-2}\xr)dx \ea\ee
owing to the following fact (with $r\ge 0$)
\be\ba\label{len4-1}
  &\int \div( |\na v|^r \na v) \div( |\na v|^2 \na v) dx\\
  &=\int \p_j(|\na v|^r \p_i v)\p_i(|\na v|^2 \p_j v)dx\\
  &=\int\p_j|\na v|^r \p_i v \p_i|\na v|^2 \p_j vdx+\int|\na v|^r \p_j\p_i v\p_i|\na v|^2 \p_j vdx\\
  &\quad+\int\p_j|\na v|^r \p_i v |\na v|^2 \p_i\p_j v dx+\int|\na v|^r \p_j\p_i v|\na v|^2 \p_i\p_j v dx\\
&=  \int \xl(2r(\na v\cdot\na^2 v)^2|\na v|^r+(r+2)|\na |\na v| |^2|\na v|^{r+2}+|\na v|^{r+2} |\na^2 v|^2 \xr)dx.
\ea\ee

Next, we have
\be\ba\label{len7}
 I_5
 &=2\ve\mu \int  |\na v|^2 \Delta v \na v \cdot u dx+ 2\ve\mu \int  \na|\na v|^2 \cdot  \na v \na v \cdot u dx\\
  &\le 32\ve\mu^2\int |\na v|^2|\Delta v|^2 dx+32\ve\mu^2\int |\na |\na v|^2|^2 dx+\frac{\ve}{16} \int   |\na v|^4|u|^2dx.\ea\ee

Notice that
\be\label{hj1} \int \n (\div u)^2dx\le 3 \int \n  |\mathcal{D} u|^2dx,\ee
this combined with H\"older inequality gives
\be\ba\label{len5}
 I_6
 \le 3\mu\int \n  |\mathcal{D} u|^2dx+\frac{\ve^2\mu}{4}\int( \div (|\na v|^2\na v) )^2 dx.\ea\ee
In order to control the last term of \eqref{len5}, we recall that  $v$ satisfies
\be\label{zqjo2}\ba  2v_t  -\ve  \div( |\na v|^2 \na v)=-  2u\cdot\na v-v\div u+\ve v^{-2p_0-1}.\ea\ee
Multiplying \eqref{zqjo2} by $\mu\ve\div (|\na v|^2\na v)  $ and integrating the resulting equality over $\O$ lead to
\be\label{zqq15}\ba &  \frac{\mu\ve}{2}\xl( \int |\na v|^4dx\xr)_t+\mu\ve^2\int(\div (|\na v|^2\na v) )^2 dx +  \mu(2p_0+1){\ve^2 } \int v^{-2p_0-2 }|\na v|^4 dx \\
&=\mu\ve  \int \div (|\na v|^2\na v) v\div u   dx+2\mu\ve  \int  (\na |\na v|^2\cdot\na v+|\na v|^2\Delta v) u\cdot\na v dx\\
&\le \frac{\ve^2\mu}{4}\int(\div (|\na v|^2\na v) )^2 dx+ \mu\int \n (\div u)^2dx+\frac{\ve}{16}\int   |\na v|^4|u|^2dx\\
&\quad+32\ve\mu^2\int   |\na v|^2|\na^2v|^2dx+32\ve\mu^2\int|\na |\na v|^2 |^2dx.\ea\ee

Submitting  \eqref{len1}, \eqref{len1'}, \eqref{len4}, \eqref{len7}, and \eqref{len5}  into \eqref{zqbb61}, then adding the resulting inequality together with \eqref{zqq15}, one has
\be\label{len6}\ba& \quad\frac{d}{dt}\xl[\frac12  \int\n |u|^2dx  +2(\kappa^2+\sqrt{\ve}\mu) \int |\na v|^2dx+\frac{\mu\ve}{2}\int |\na v|^4dx  \xr]  \\
&\quad+2\nu\int \n  |\mathcal{D}u|^2  dx +\frac{\sqrt{\ve}}{2}\int \n  |\na u|^2  dx +\frac{\ve}{8}\int   |\na v|^4 |u|^2dx
+\frac{\ve}{2}\int \n^{-p_0} |u|^2dx\\
&\quad+2(\kappa^2+\sqrt{\ve}\mu)\ve\int \xl(|\na v|^2 |\na^2 v|^2+\frac{1}{2}|\na |\na v|^2|^2+(2p_0+1)|\na v|^2v^{-2p_0-2}\xr)dx\\
&\quad+\frac{1}{2}\mu\ve^2\int( \div (|\na v|^2\na v) )^2 dx +  \mu(2p_0+1){\ve^2 } \int v^{-2p_0-2 }|\na v|^4dx\\
&\quad+\ve^{3/2}\int\n  |w|^{3}|u|^{2}dx +r_0\int |u|^{2}dx+r_1\int\n  |u|^{4}dx+\int   u \cdot \na\n^\ga  dx+\ve\mu\int   u \cdot \na\n^{-p_0}  dx\\
&\quad\le 6\mu\int \n  |\mathcal{D} u|^2dx+ 128\ve \mu^2\int |\na|\na v|^2|^2 dx + 88\ve \mu^2 \int |\na^2 v|^2|\na v|^2 dx.\ea\ee

Now, for the last two  terms on the left hand side  of \eqref{len6},   it holds that for $q\not=1,$
\be\label{zqbini5}\ba    \int  u\cdot\na  \n^{q }  dx    &= -\frac{q}{q -1}\int\n^{q -1}\div(\n u) dx\\ &= -\frac{q}{q -1}\int\n^{q -1}(-\n_t  +\ve v\div(|\na v|^2\na v)+\ve \n^{-p_0}) dx \\ &= \frac{1}{   q-1 }\xl(\int\n^{q }dx\xr)_t+\frac{q(2q -1)\ve}{ q -1 }\int\n^{q -1} |\na v|^4dx -\frac{q\ve}{q -1}\int\n^{q -1-p_0}dx.\ea\ee
Choosing $q=-p_0 $ in \eqref{zqbini5}, one gets
\be \label{zqjo1}\ba     &\frac{1}{6(  p_0+1) }\xl(\int\n^{-p_0 }dx\xr)_t+\frac{  p_0(2p_0+1)\ve^2 }{6(  p_0+1) }\int\n^{-p_0 -1} |\na v|^4dx +\frac{  p_0\ve^2 }{6(  p_0+1)}\int\n^{ -1-2p_0}dx \\&= \frac{\ve }{6} \int  \n^{-p_0 } \div u  dx\\&\le  \frac{p_0\ve^2 }{12(p_0+1)}\int\n^{ -1-2p_0}dx+\frac{1}{2}\int \n |\mathcal{D} u|^2dx .\ea\ee

 Finally, choosing
\be\label{amu2}  11\kappa \le  \nu  \ee
such that
\be\label{amu1} 20\mu< \nu, ~~~
400\mu^2<  \kappa^2,  \ee
multiplying \eqref{zqjo1} by $\nu$ and $6\mu$, respectively, then adding the resulting inequalities, \eqref{zqbc6} and \eqref{len6}  together,
 we thus obtain \eqref{zqbb30} after  using \eqref{zqbini5}, \eqref{len4-1}, \eqref{amu1},  Gronwall's inequality, and the following simple fact  \bnn \n^{-p_0+\ga-1}\le \n +\n^{-p_0}.\enn
Hence, the proof of Lemma \ref{zqlem10} is finished.  \hfill$\Box$

Next, with the same spirit of the BD  entropy estimates   due to Bresch-Desjardins \cite{BreschDesjardinsLin05,BreschDesjardins03,
BreschDesjardins03b,BreschDesjardins02}, we have the following estimates in Lemma \ref{zqlem11}.

\begin{lemma}
\label{zqlem11} There  exists some generic constant $C$ independent of $\ve$, $r_0$, $r_1$, and $\kappa$ such that
\be\label{zqbb5'}\ba& \sup_{0\le t\le T}\int  \left( |\na v|^2+\ve |\na v|^4-r_0\log_{-}\n   \right)dx   + \int_0^T \int \left( \n   |\na u|^2 +|\na(\sqrt{\n}u)-u\otimes \na\sqrt{\n}|^2 +   \n^{\ga-2}   |\na \n|^2\right)dxdt\\
 &+(\ka^2+\sqrt{\ve}\mu)\int_0^T \int \n |\na^2\log \n|^2 dxdt+\ve\nu\int_0^T\int \xl(|\na v|^2 |\na^2v|^2+|\na v|^2 |\na |\na v||^2  +\n^{-p_0-1}|\na v|^2\xr)dxdt\\&+\ve^2 \int_0^T\int \xl(|\na v|^4|\na^2v|^2+|\na v|^4 |\na |\na v||^2+\n^{-p_0-1}|\na v|^4\xr)dxdt+\ve\mu\int_0^T \int v^{-2}|\na v|^6dxdt\\
 &+r_0\ve\int_0^T\int \xl(v^{-2} |\na v|^4 +\n^{-p_0-1}\xr)dxdt \le C+Cr_0 +Cr_1.\ea\ee
 Furthermore, it holds that
 \begin{align}\label{useful}
 \ve^{\frac{3}{2}} \int_0^T\int \xl(\rho |w|^5+ \rho |u|^5\xr) dxdt + \ve   \int_0^T\int  \xl(v^{-2 } |\nabla v|^6 + v^{-3}|\nabla v|^5 \xr) dxdt \le C+Cr_0 +Cr_1.
 \end{align}
\end{lemma}
{\it Proof. } First,  set
\be\notag \label{zqlalm1} G  \triangleq \ve v\div( |\na v|^2 \na v)+\ve\n^{-p_0},\ee
multiplying \eqref{zqba1}$_1$ by $\n^{-1}$ and applying gradient to the resulting equality lead to
\be  \label{hj2} (\na \log\n)_t+u\cdot\na\na\log\n+\na u\cdot\na\log\n+\na \div u=\na(\n^{-1}G).\ee
Thus, multiplying \eqref{hj2} by $\na \n$, we obtain after using integration by parts and \eqref{zqba1}$_1$ that
\be \ba \label{hj3} \frac{1}{2} \xl(\int \n^{-1}|\na \n|^2dx\xr)_t&+\int\n^{-1}\na\n\cdot\na u\cdot
\na \n dx+\int\na\n\cdot\na\div udx\\&+\int \n^{-1}G\xl(\Delta\n-\frac{1}{2}\n^{-1}|\na \n|^2\xr)dx=0.\ea\ee
 Then, multiplying \eqref{zqba1}$_2$ by $\na\log\n=\n^{-1}\na \n$ and integrating by parts yield
 \be \ba \label{hj4}  &\int u_t \xd \na\n dx+\int u\xd \na u\xd \na \n dx-2\nu\int \div (\n \mathcal{D}u) \xd \na\log\n  dx-\sqrt{\ve}\int \div(\n\na u)\xd \na\log\n dx \\
 &+\int P'(\n)\n^{-1}|\na \n|^2dx+ 2(\ka^2+\sqrt{\ve}\mu)\int \n |\na^2\log \n|^2 dx\\
 &=  \ve\int v |\na v|^2 \na v\cdot \na u\cdot\na \log \n  dx+ \ve\mu\int v |\na v|^2 \na v\cdot \na (\na \log\n)\cdot\na \log \n  dx\\
&\quad- \ve\int  \n^{-p_0} u\cdot\na \log \n  dx- \ve^{3/2}\int \n  |w|^{3}u  \cdot\na \log \n  dx \\
 &\quad-r_0\int u\cdot\na\log\n dx-r_1\int\n|u|^2u\cdot\na\log\n dx  - \ve\mu\int\na( v\div( |\na v|^2 \na v))\cdot\na \log \n dx\\
 &\quad + \ve\mu\int v\div( |\na v|^2 \na v)\na \log \n \cdot\na \log \n dx- \ve\mu\int  \na \n^{-p_0} \cdot\na \log \n dx\\
 &\triangleq \sum_{i=1}^{9} \tilde{I}_i,\ea\ee
where the first term on the left hand of \eqref{hj4} can be handled as follows
\be \ba \label{hj5}  \int u_t \xd \na\n dx&=\xl(\int u\cdot \na \n dx\xr)_t-\int u\cdot\na u\cdot\na \n dx\\
&\quad-2\int \n \mathcal{D}u:\na udx+\int \n |\na u|^2dx+\int \div u Gdx.\ea\ee

Adding \eqref{hj3} multiplied by $2\nu+\sqrt{\ve}$ to \eqref{hj4} and using \eqref{hj5}, one has
\be\label{zqbb7}\ba& \frac{2\nu+\sqrt{\ve} }{2}\xl(\int\n^{-1} |\na  \n |^2dx\xr)_t + \xl(\int   u\cdot \na  \n  dx \xr)_t+ \int  \n |\na u|^2dx+2(\ka^2+\sqrt{\ve}\mu)\int \n |\na^2\log \n|^2 dx\\&  \quad + \int P'(\n) \n^{-1}|\na \n|^2dx +(2\nu+\sqrt{\ve} )\int  \n^{-1} G\left( \Delta   \n -\frac12  \n^{-1}|\na  \n|^2 \right)dx \\
&=- \int  G\div u dx+2  \int  \n \mathcal{D}u: \na  u   dx +\sum_{i=1}^{9} \tilde{I}_i.\ea\ee

Since \be  \label{hi5}\Delta   \n -\frac12  \n^{-1}|\na  \n|^2 =2v\Delta v,\ee the last term on the left-hand side of \eqref{zqbb7} can be calculated  as
\be\label{qq3}\ba &\int\ \n^{-1} G\left( \Delta   \n -\frac12  \n^{-1}|\na  \n|^2 \right)dx\\&=2\ve\int   \div (|\na v|^2\na v) \Delta v dx+2\ve\int \n^{-p_0-1/2}\Delta v dx \\&= 2\ve\int    |\na v|^2|\na^2v|^2dx+ \ve\int  |\na|\na v|^2|^2dx  +2(2p_0+1)\ve\int \n^{-p_0-1 }|\na v|^2 dx ,\ea\ee
where we have used \eqref{len4-1} with $r=0$. 

 Now, we will estimate each term on the righthand side of \eqref{zqbb7} in the following way.

 First, with the same arguments as those in \cite{li04}, one has
\be \label{zqq5}\ba  &-\int \div u Gdx+2\int \n \mathcal{D}u:\na udx+\tilde{I}_1+\tilde{I}_3\\
&\le  \frac{\ve^2}{8}\int(\div (|\na v|^2\na v))^2dx+ \frac{\nu}{2}\ve \int|\na v|^2|\na^2v|^2dx+\frac{1}{4 }\int   \n | \na  u|^2   dx \\
&\quad+C\ve^2\int \n^{-2p_0-1}dx+C(\nu)\ve \int  |\na v|^4| u|^2dx+C\int  \n  |\mathcal{D} u|^2 dx.\ea\ee

Next, it holds
\be \ba \label{zzqq14} \tilde{I}_{2}&=  -\frac{1}{2}\ve\mu\int \div(v |\na v|^2 \na v)|\na \log \n|^2  dx\\
&=-2\ve\mu\int v^{-2}|\na v|^2\xl(|\na v|^4+ v\na |\na v|^2\cdot\na v+v|\na v|^2\Delta v\xr) dx\\
&=-2\ve\mu\int v^{-2}|\na v|^6dx -2\ve\mu\int v^{-1}|\na v|^2\na |\na v|^2\cdot\na vdx-2\ve\mu\int v^{-1}|\na v|^4\Delta vdx\\
 &\le -\ve\mu\int v^{-2}|\na v|^6dx+2\ve\mu\int |\na |\na v|^2|^2dx+2\ve\mu\int  |\na v|^2|\Delta v|^2dx\ea\ee




Recalling the definition of $w$   and using Young's inequality, one gets
\be \ba\label{lbd1}  \tilde{I}_4
&=-2\ve^{3/2}\int \n^{1/2}|w|^3u\cdot\na vdx \le   \ve^{3/2}\int \n |w|^3|u|^2dx
+\ve^{3/2}\int  |w|^3|\na v|^2dx\\
&\le  C(\nu)\ve^{3/2}\int \n |w|^3|u|^2dx
+ C(\nu)\ve^{3/2}\int  v^{-3}|\na v|^5dx,
\ea\ee
where in the last inequality we have used the following fact:
\be \ba\label{lbd2}
 &\ve^{3/2}\int   |w|^3|\na v|^2dx\\
 &\le\frac{1}{16\nu^2} \ve^{3/2}\int \n|w|^3|u+\mu\na \log\n|^2dx+C(\nu)\ve^{3/2}\int \n^{-3/2}|\na v|^5dx\\
&\le \frac{1}{8\nu^2}\ve^{3/2}\int \n|w|^3|u|^2dx+ \frac{1}{2}\ve^{3/2}\int|w|^3|\na v|^2dx+C(\nu)\ve^{3/2}\int v^{-3}|\na v|^5dx.\ea\ee
 The last term on the left hand of \eqref{lbd1} can be handled as follows:
\be \ba\label{lbd5}    \int v^{-3}|\na v|^5dx&=  \int v^{-3}|\na v|^3\na v\cdot\na vdx\\
&=3\int v^{-3}|\na v|^5dx- \int v^{-2}\xl(\na|\na v|^3\cdot\na v+|\na v|^3\Delta v\xr)dx,\ea\ee
which along with \eqref{len4-1} and Young's inequality shows
\be \ba\label{lbd6}    C(\nu)\ve^{3/2}\int v^{-3}|\na v|^5dx& =C(\nu)\ve^{3/2}\int v^{-2}\xl(\na|\na v|^3\cdot\na v+|\na v|^3\Delta v\xr)dx\\
&\le   \frac{1}{8}\ve^2 \int \xl(|\na v|^4|\na|\na v| |^2+|\na v|^4|\na^2 v|^2\xr)dx+C(\nu)\ve \int \n^{-2} |\na v|^2dx\\
&\le \frac{\ve^2}{8}\int(\div (|\na v|^2\na v))^2dx+\ve \int |\na v|^4dx\\
&\quad+C(\nu)\ve  \int \n dx+C(\nu)\ve  \int \n^{-p_0} dx.\ea\ee
Combined this with \eqref{lbd1} yields that \be \ba\label{lbd7}  \tilde{I}_4&\le  C(\nu)\ve^{3/2}\int \n |w|^3|u|^2dx+\frac{\ve^2}{8}\int(\div (|\na v|^2\na v))^2dx\\
&\quad+C(\nu) \int \n dx+C(\nu)\ve  \int \n^{-p_0} dx+\ve \int |\na v|^4dx.
\ea\ee

The terms $\tilde{I}_5$--$\tilde{I}_{8}$ can be handled by some directly calculations:
\be \ba\label{lbd8}  \tilde{I}_5 &=-r_0\int \frac{u\cdot\na \n}{\n} dx= r_0\int \frac{\n_t+\n\div u -\ve v\div( |\na v|^2 \na v)-\ve\n^{-p_0}}{\n} dx \\
&=r_0\xl(\int \log\n dx\xr)_t  -r_0\ve\int v^{-2}|\na v|^4dx-r_0\ve\int \n^{-p_0-1}dx,
\ea\ee
\be \ba\label{gylbd8}   \tilde{I}_6 
 = r_1\int |u|^2\div u\n dx+2r_1\int u\cdot\na u\cdot u  \n dx\le  Cr_1^2\int \n|u|^4dx+ \frac{1}{4}\int\n|\na u|^2 dx,
\ea\ee
and
\be \ba\label{lbd13}  \tilde{I}_{7} +\tilde{I}_{8} +\tilde{I}_{9} 
&\le \frac{\ve^2}{4}\int |\div( |\na v|^2 \na v)|^2dx+\frac{3\mu^2}{2}\int \n |\na^2 \log \n|^2 dx+\frac{\ve^2}{2}\int \n^{-2p_0-1}   dx \\
&\quad+\frac{\ve\mu}{2}\int v^{-2}|\na v|^6dx +8\ve\mu\int |\na |\na v|^2|^2dx+8\ve\mu\int |\na^2 v|^2|\na v|^2dx.\ea\ee

Substituting \eqref{qq3}--\eqref{zzqq14} and \eqref{lbd7}--\eqref{lbd13} into \eqref{zqbb7}, we obtain after using \eqref{amu2}--\eqref{amu1}  that
\be\label{lbd001}\ba& \frac{2\nu+\sqrt{\ve} }{2}\xl(\int\n^{-1} |\na  \n |^2dx\xr)_t + \xl(\int   u\cdot \na  \n  dx \xr)_t  +\frac{1}{2}\int   \n | \na  u|^2   dx\\
&  \quad + \xl(\frac{\ka^2}{2}+\sqrt{\ve}\mu\xr)\int \n |\na^2\log \n|^2 dx+ \int P'(\n) \n^{-1}|\na \n|^2dx +\frac{1}{2}\ve\mu\int v^{-2}|\na v|^6dx \\
& \quad+(2\nu+2\sqrt{\ve} )\ve\xl( \int    |\na v|^2|\na^2v|^2dx+ \frac{1}{2}\int  |\na|\na v|^2|^2dx +(2p_0+1)\int \n^{-p_0-1 }|\na v|^2 dx \xr)  \\
&\quad+r_0\ve\int v^{-2}|\na v|^4dx+r_0\ve\int \n^{-p_0-1}dx\\
&\le \frac{\ve^2}{2} \int( \div(|\na v|^2\na v))^2dx + C\ve^2 \int \n^{-2p_0-1} dx +C\int \n \xl(|\mathcal{D}u|^2+|\div u|^2\xr)dx\\
&\quad+\frac{\nu}{2}\ve\int|\na v|^2|\na^2v|^2dx+C(\nu)\ve \int  |\na v|^4| u|^2dx+C(\nu)\ve^{3/2}\int \n |w|^3|u|^2dx \\
&\quad+C(\nu) \int \n dx+C(\nu)\ve  \int \n^{-p_0} dx
+\ve \int |\na v|^4dx+Cr_1^2\int \n|u|^4dx+r_0\xl(\int \log\n dx\xr)_t.
\ea\ee

 Next, with the similar arguments as \eqref{zqq15}, it holds that
\be\label{qq152}\ba & \xl( \frac12\int \ve|\na v|^4dx\xr)_t+\ve^2\int( \div (|\na v|^2\na v) )^2 dx  +  (2p_0+1){\ve^2 } \int v^{-2p_0-2 } |\na v|^4 dx \\& =  {\ve } \int  \div (|\na v|^2\na v)   v\div u   dx+2{\ve } \int \div (|\na v|^2\na v) u\cdot\na v dx\\
&\le \frac{\ve^2}{4}\int( \div (|\na v|^2\na v) )^2 dx+C\int \n (\div u)^2dx + \frac{\nu\ve}{2} \int|\na v|^2|\na^2v|^2dx +C(\nu) \ve\int |u|^2|\na v|^4dx .\ea\ee

 The combination of \eqref{lbd001} with \eqref{qq152} yields
 \be\label{lbd002}\ba& \frac{2\nu+\sqrt{\ve} }{2}\xl(\int\n^{-1} |\na  \n |^2dx\xr)_t + \xl(\int   u\cdot \na  \n  dx \xr)_t+ \frac\ve2\xl( \int |\na v|^4dx\xr)_t+\frac{1}{2}\int   \n | \na  u|^2   dx\\
 &\quad+ \frac{\ve^2}{4}\int( \div (|\na v|^2\na v) )^2 dx  +  (2p_0+1){\ve^2 } \int v^{-2p_0-2 } |\na v|^4 dx\\
&  \quad + \xl(\frac{\ka^2}{2}+\sqrt{\ve}\mu\xr)\int \n |\na^2\log \n|^2 dx+ \int P'(\n) \n^{-1}|\na \n|^2dx +\frac{1}{2}\ve\mu\int v^{-2}|\na v|^6dx \\
& \quad+  ( \nu+2\sqrt{\ve} )\ve\xl( \int    |\na v|^2|\na^2v|^2dx+ \frac{1}{2}\int  |\na|\na v|^2|^2dx +(2p_0+1)\int \n^{-p_0-1 }|\na v|^2 dx \xr)  \\
&\quad+r_0\ve\int v^{-2}|\na v|^4dx+r_0\ve\int \n^{-p_0-1}dx\\
&\le   C\ve^2 \int \n^{-2p_0-1} dx +C\int \n |\mathcal{D}u|^2 dx +C(\nu)\ve \int  |\na v|^4| u|^2dx+C(\nu)\ve^{3/2}\int \n |w|^3|u|^2dx\\
&\quad +C(\nu) \int \n dx+C(\nu)\ve  \int \n^{-p_0} dx +\ve \int |\na v|^4dx+Cr_1^2\int \n|u|^4dx+r_0\xl(\int \log\n dx\xr)_t
  \\
&\triangleq H+r_0\xl(\int \log\n dx\xr)_t,\ea\ee
 On the one hand, one deduces from \eqref{zqbb30} that  $H$  satisfies
\be\label{lbd003}\ba \int_0^T Hdt\le C +Cr_1. \ea\ee
On the other hand, recalling that $-  \log_{-}\n_0 \in L^1$ in \eqref{pini1}$_3$ and using \eqref{zqbb30}, it holds
\be\label{gylbd003}\ba \int_0^T r_0\xl(\int \log\n dx\xr)_tdt&=r_0\int \log\n dx -r_0\int \log\n_0 dx\\
&=r_0\int \log_{-}\n dx +r_0\int \log_{+}\n dx-r_0\int \log_{+}\n_0 dx-r_0\int \log_{+}\n_0 dx\\
&\le r_0\int \log_{-}\n dx    +Cr_0,\ea\ee
where $\log_{+}g\triangleq\log \max\{1,g\}$.

Noting that
\be\label{gyg3}\sqrt{\n} \na u=\na(\sqrt{\n}u)-u\otimes  \na \sqrt{\n},\ee
we thus deduce \eqref{zqbb5'} directly by  integrating \eqref{lbd002} over $[0,T]$ and   using  \eqref{lbd003},  \eqref{gylbd003}, \eqref{zqbb30}, and \eqref{gyg3}.

Finally, some directly calculations together with H\"older inequality and \eqref{len1-5} deduce that
\begin{align*}
& \ve   \int   v^{-3}|\nabla v|^5 dxdt +\ve\int v^{-2}|\na v|^6dx \\
&\le C\ve \int v^{-2}|\na v|^6dx+C \ve\int \n^{-4}dx\\
&\le C\ve \int|\na v|^2 |\Delta v|^2dx +C\ve\int|\na |\na v|^2 |^2dx+C \ve  \int\xl(\n+\n^{-p_0}\xr)dx,
\end{align*}
which along with  (\ref{zqbb30}) and (\ref{zqbb5'}) shows that
\be\label{xims2} \ba
\ve \int_0^T\int  v^{-2 } |\nabla v|^6 dxdt +\ve  \int_0^T\int   v^{-3}|\nabla v|^5 dxdt\le C+Cr_0+Cr_1.
\ea\ee
Then it follows from \eqref{zqbb30}, (\ref{zqbb5'}),  \eqref{xims2}, and H\"older inequality that
\be\label{lowb2}\ba  &\ve^{3/2}\int_0^T\int (\n |w|^{5}+\n |u|^{5}) dx \\
&= \ve^{3/2} \int_0^T\int \xl(\n |w|^3|u+\mu\na\log\n|^2+\n |w-\mu\na\log\n|^3|u|^2 \xr)dx \\
&\le C\ve^{3/2} \int_0^T\int \n |w|^3|u|^2 dx +C\ve^{3/2} \int_0^T\int \n  |w|^3|\na\log\n|^2dxdt+ C \ve^{3/2} \int_0^T\int \n |\na\log\n|^3|u|^2 dx \\
&\le C+Cr_0+Cr_1+ C \ve^{3/2} \int_0^T \int  |\na v|^5 v^{-3}dx+\frac{1}{2}\ve^{3/2} \int_0^T\int \n |w|^{5} dx+\frac{1}{2}\ve^{3/2} \int_0^T\int \n |u|^{5} dx\\
&\le C+Cr_0+Cr_1+\frac{1}{2}\ve^{3/2} \int_0^T\int \n |w|^{5} dx+\frac{1}{2}\ve^{3/2} \int_0^T \int \n |u|^{5} dx.\ea\ee

Thus, The combination of \eqref{xims2} and \eqref{lowb2} gives (\ref{useful}). The proof of Lemma \ref{zqlem11} is completed.  \hfill $\Box$

Now, using the BD-entropy inequality obtained in Lemma \ref{zqlem11},
we can obtain following useful a priori estimates.

\begin{lemma}
\label{newlem} There  exists some generic constant $C$ independent of $\ve$,  $r_0$, $r_1$, and $\kappa$ such that
\be\label{newb}\ba&  \kappa^2\int_0^T\int\xl(|\nabla\rho^\frac14|^4+|\nabla^2\rho^\frac12|^2\xr)dxdt
+r_1\kappa\int_0^T \int  |\na (\sqrt{\n} u )|^2  dxdt \le C+Cr_0+Cr_1. \ea\ee
\end{lemma}
{\it Proof.} First, recalling the following facts due to J\"ungel \cite{jun} (see also \cite[Lemma 2.1]{vy1})
\be\label{newb1} \ba
  \int|\nabla\rho^\frac14|^4dx\le 8 \int\rho|\nabla^2\log\rho|^2dx, ~~~
  \int|\nabla^2\rho^\frac12|^2dx\le 7\int\rho|\nabla^2\log\rho|^2dx,  \ea\ee
which combined with  \eqref{zqbb5'} shows that
\be\label{gynewb1} \ba
 \kappa^2\int_0^T\int\xl(|\nabla\rho^\frac14|^4+|\nabla^2\rho^\frac12|^2\xr)dxdt
\le
15 \kappa^2\int_0^T\int\rho|\nabla^2\log\rho|^2dxdt\le C+Cr_0+Cr_1.   \ea\ee

Next, we have
\begin{align}\label{1}
\nabla(\sqrt{\rho}u)  =\sqrt{\rho}\nabla u+u\otimes \nabla\sqrt{\rho}
& =\sqrt{\rho}\nabla u+2\rho^{\frac14}u\otimes  \nabla\rho^{\frac14},
\end{align}
which along with \eqref{gynewb1}, \eqref{zqbb30}, and \eqref{zqbb5'} that
\be\ba \label{newb2}\notag
\int_0^T \int r_1\kappa|\na (\sqrt{\n} u )|^2  dxdt&\le 2r_1\kappa\int_0^T \int   \rho |\nabla u|^2dxdt+8r_1\kappa\int_0^T \int  \rho^{\frac12} |u|^2|\nabla\rho^{\frac14}|^2dxdt\\
&\le Cr_1 \int_0^T \int   \rho |\nabla u|^2dxdt+4r_1^2\int_0^T \int \rho |u|^4dxdt+4\kappa^2\int_0^T \int|\nabla\rho^{\frac14}|^4dxdt\\
&\le C+Cr_0+Cr_1.
\ea\ee
This combined with \eqref{gynewb1} gives \eqref{newb} and thus finishes the  proof of Lemma \ref{newlem}.  \hfill$\Box$

Following the same arguments as those in \cite{li04}, we will use  a De Giorgi-type procedure to obtain the following estimates on the lower and upper  bounds  of the density which are crucial   to obtain the global existence of strong solutions to the problem \eqref{zqba1}--\eqref{zqbb35}.

\begin{lemma}
\label{zqlem13}There exists some positive constant $C $ depending on $\ve$, $r_0$,  $r_1$, and $\kappa$ such that for all $(x,t)\in \O\times (0,T)$ \be \label{qc1} C^{-1}\le \n(x,t)\le C . \ee \end{lemma}
{\it Proof.} The proofs are similar to the arguments in Li-Xin \cite[Lemma 4.4]{li04}. We sketch them here for completeness.

First, it follows from \eqref{zqbb5'}, \eqref{zqbb30}, and  Sobolev inequality   that
\be\label{qq6}\ba \sup_{0\le t\le T}\|\n\|_{L^\infty}&=\sup_{0\le t\le T}\|  v\|^2_{L^\infty}\le C\sup_{0\le t\le T}\left(\|  v\|_{L^{2 }}+\|\na v\|_{L^4}\right)^2\le \hat C. \ea\ee

Next, we will use a De Giorgi-type procedure to obtain the     lower   bound  of the density. In fact,   since $h\triangleq v^{-1}$ satisfies
\be \label{pbb81} \ba&2h_t+ 2u\cdot \na h- h\div u+\ve h^{2p_0+3} +2\ve h^{-5}|\na h|^4 =  \ve \div (h^{-4}|\na h|^2\na h),\ea\ee
multiplying \eqref{pbb81} by $(h-k)_+$ with $ k\ge \|h(\cdot,0)\|_{L^\infty } =\|\n_0^{-1/2}\|_{L^\infty } $ yields that
\be\label{pbb82}\ba & \sup_{0\le t\le T}\int  (h-k)_+^2dx
+\ve\int_0^T\int h^{-4}|\na(h-k)_+ |^4 dxdt \\
&\le C\int_0^T\int h|u||\na (h-k)_+|dxdt+C\int_0^T\int (h-k)_+|u||\na h|dxdt\\
&\le C\int_0^T\int 1_{\hat A_k} \n^{-4/3}|u|^{ 4/3} dxdt
+\frac{\ve}{2} \int_0^T\int h^{-4} |\na(h-k)_+ |^4dxdt, \ea\ee
where $\hat A_k\triangleq \{(x,t)\in \O\times (0,T)|h(x,t)>k\}.$
Denote $\hat \nu_k\triangleq  |\hat A_k|$,
it follows from H\"older  inequality, \eqref{zqbb30}, and \eqref{useful}  that
\be\label{pbb83}\ba &\int_0^T\int 1_{\hat A_k} \n^{-4/3}|u|^{ 4/3} dxdt  \\ &\le \left(\int_0^T\int 1_{\hat A_k}\n^{-24/11}  dxdt\right)^{11/15}\left(\int_0^T\int \n |u|^{5} dxdt\right)^{4/15} \\&\le C\left(\int_0^T\int(\n+\n^{-p_0})  dxdt\right)^{1/15}| \hat A_k |^{2/3 }\\&\le C\hat \nu_k^{2/3}.\ea\ee

Now,  submitting \eqref{pbb83} into \eqref{pbb82}   leads to
\be\label{gai1}\ba \sup_{0\le t\le T}\int (h-k)_+^2dx+  \int_0^T\int  h^{-4}|\na(h-k)_+ |^4dxdt \le C\hat \nu_k^{2/3}.\ea\ee
This together with \eqref{zqbb30} and H\"older  inequality gives
\be\ba\label{gai2}
\int_0^T\int|\nabla(h-k)_+|^2dxdt
&= \int_0^T\int 1_{\hat{A}_k}h^2h^{-2}|\nabla(h-k)_+|^2dxdt \\
& \leq \int_0^T \xl(\int 1_{\hat{A}_k}^{3}dx\xr)^{1/3}\xl(\int  h^{12}dx\xr)^{1/6}
\left(\int  h^{-4}|\nabla(h-k)_+|^4dx\right)^{1/2}dt \\
& \leq C\hat{\nu}_k^{1/3}
 \xl(\int_0^T \left(\int  h^{12}dx\right)^{1/3}dt\xr)^{1/2}
\left(\int_0^T\int  h^{-4}|\nabla(h-k)_+|^4dxdt\right)^{1/2} \\
& \le C\hat{\nu}_k^{2/3}\xl(\int_0^T \left(\int  (\n+\n^{-p_0})dx\right)^{1/3}dt\xr)^{1/2}\\
& \le C\hat{\nu}_k^{2/3}.
\ea\ee

Hence, the  Sobolev inequality combined with \eqref{gai1} and \eqref{gai2} derive that
\be\label{pbb84}\ba \|(h-k)_+\|_{L^{10/3}(\O\times (0,T))}^2&\le C   \sup_{0\le t\le T}\int (h-k)_+^2dx+ C \int_0^T\int  |\na(h-k)_+ |^2dxdt \le C \hat \nu_k^{2/3 }.\ea\ee
This implies that for $\bar{k}>k,$   \be \label{qq8}\ba \hat \nu_{\bar{k}} \le C(\bar{k}-k)^{-10/3}\hat\nu_{k}^{10/9 }\ea\ee due to the following simple fact that
\bnn (\bar{k}-k)^2 |\hat A_{\bar{k}}|^{3/5}\le \|(h-k)_+\|_{L^{10/3}(\O\times (0,T))}^2.\enn

Finally, it follows from \eqref{qq8} and the  De Giorgi-type lemma \cite[Lemma 4.1.1]{wyw} that there exists some positive constant $C\ge \hat C$ such that   \bnn \sup_{(x,t)\in \O\times (0,T)}\n^{-1}(x,t)\le   C,\enn
which along with \eqref{qq6} gives \eqref{qc1} and thus completes  the proof of Lemma \ref{zqlem13}.  \hfill$\Box$

In order to overcome the difficulties come  from the third order tensor term in \eqref{qns}$_2$, we will use a transformation through the effective velocity $w$ which is defined in \eqref{lw}.
 Next lemma shows that the
system of $(\n,u)$ can be  written equivalently in terms of $(\n,w)$.

\begin{lemma}\label{uw0} Let $(\n,u)$ be a smooth solution of the system \eqref{zqba1}, then $(\n,w)$ with $w$ defined in \eqref{lw} will satisfy the following system
\be\label{w=u0} \begin{cases}   \n_t+\div(\n w)    =\mu \Delta \n+  \ve v\div( |\na v|^2 \na v)+\ve\n^{-p_0},
     \\
     \n w_t+  \n w\cdot\na w   +\na P  -2(\nu-\mu)\div( \n    \mathcal{D} w) - \mu\n \Delta w- \sqrt{\ve}\div(\n \na w)  \\
  \quad= 2\mu\na\n \cdot\na w +\ve v|\na v|^2\na v\cdot\na w- \ve^{3/2}\n    |w|^{3}u  -r_0u-r_1\n|u|^2u  -\ve\n^{-p_0} w.   \end{cases}\ee
\end{lemma}
{\it Proof.}
First, it is easy to deduce from \eqref{zqba1}$_1$ that
\be \label{w=u1}\n_t+\div(\n w)    =\mu  \Delta \n+  \ve v\div( |\na v|^2 \na v)+\ve\n^{-p_0}.\ee

In order to prove \eqref{w=u0}$_2$,   we recall some   identities as follows:
\be \begin{cases}  \label{w=u2}
\mu(\n\na \log\n)_t=-\mu\na\div(\n u)+ \ve \mu\na (v\div( |\na v|^2 \na v))+\ve\mu \na\n^{-p_0},\\
\mu\div(\n u\otimes\na\log\n+\n\na\log\n\otimes u)=\mu\Delta(\n u)-2\mu \div(\n\mathcal{D}u)+\mu\na\div(\n u),\\
\mu^2 \div( \n\na\log\n\otimes \na\log\n)=\mu^2\Delta(\n \na \log \n)-\mu^2\div(\n\na^2\log\n).
\end{cases} \ee
Fuethermore,  using \eqref{lw} and \eqref{len0},  one can  rewrite  \eqref{zqba1}$_2$  as
\be \label{w=u3} \ba &(\n u)_t+ \div(\n u\otimes u)   -2\nu\div( \n    \mathcal{D} u) +\na P +\ve\mu \na\n^{-p_0}- \ka^2\div(\n \na^2\log \n)+r_0u+r_1\n|u|^2u\\
&=  \sqrt{\ve}\div(\n \na w)  - \ve\mu\na( v\div( |\na v|^2 \na v))+\ve v|\na v|^2\na v\cdot\na w
+\ve v\div( |\na v|^2 \na v) w - \ve^{3/2}\n |w|^{3}u.
\ea\ee
Notice that $\mu=\nu-\sqrt{\nu^2-\kappa^2}$, one thus obtains after  adding  \eqref{w=u2} and \eqref{w=u3} together that
\be \label{w=u4} \ba &(\n w)_t+ \div(\n w\otimes w)  +\na P  -2(\nu-\mu)\div( \n    \mathcal{D} w)-\mu\Delta(\n w)- \sqrt{\ve}\div(\n \na w)  \\
&= (\ka^2-\mu^2-2(\nu-\mu)\mu)\div(\n \na^2\log \n) +\ve v|\na v|^2\na v\cdot\na w+\ve  v\div( |\na v|^2 \na v) w\\
&\quad-r_0u-r_1\n|u|^2u- \ve^{3/2}\n    |w|^{3}u  \\
&= \ve v|\na v|^2\na v\cdot\na w+\ve  v\div( |\na v|^2 \na v) w -r_0u-r_1\n|u|^2u- \ve^{3/2}\n    |w|^{3}u .\ea\ee
This combined with \eqref{w=u1} gives directly \eqref{w=u0} and finishes the proof of Lemma \ref{uw0}.  \hfill$\Box$

 Next, with the estimates of $(\n, u)$ in Lemmas \ref{zqlem10}--\ref{zqlem13} in hand, we will derive  some estimates on $(\n, w)$ in following Lemma \ref{hiu}.

\begin{lemma}\label{hiu} There exists some constant $C $ depending on $\ve$, $r_0$,  $r_1$, and $\kappa$ such that \be\label{hi0}\ba  &\sup_{0\le t\le T}  (\|w\|_{L^2\cap L^4}+\|\na v\|_{L^2\cap L^4 }) \\
&+\int_0^T\int \left( |w |^{5}+|w |^{7}+|\na v|^4|\na^2v|^2 + |\na w |^2+|\na v|^{21}\xr)dxdt\\
&+\int_0^T\int \left( |\na w|^2|w|^2 +|\div w|^2|w|^2 +|u|^4|w|^2+|u|^2|w|^4\right)dxdt\le C. \ea\ee
\end{lemma}
{\it Proof.} First, it follows from \eqref{qc1},   \eqref{zqbb30},  \eqref{zqbb5'},  and \eqref{useful} that \be \label{zqbb89}\sup_{0\le t\le T}  (\|w\|_{L^2  }+\|\na v\|_{L^2\cap L^4 }) +\int_0^T\int \left( |\na v|^4|\na^2v|^2 +|\na v|^6+|\na w |^2+ |w |^{5}\right)dxdt\le C. \ee

Then it   follows from \eqref{w=u0}$_1$ and \eqref{hi5} that  $v$ satisfies
   \be \label{hi6} \ba
2 v_t-2\mu\Delta v-\ve\div(|\na v|^2\na v)  =-v\div w-2w\cdot\na v+2\mu v^{-1}|\na v|^2+\ve v^{-2p_0-1}.
\ea\ee
This yields that
\be\label{zqbb87}2 v_t -\ve \div((2\mu\ve^{-1}+|\na v|^2)\na v) =-\div(w  v+\na g)- \frac{1}{|\O|}\int ( w\cdot \na v-\ve v^{-2p_0-1}-2\mu v^{-1}|\na v|^2) dx,  \ee
where   $g(\cdot,t)$ (with $t>0$)  is the unique solution to the following problem
\be\label{qbb88}\begin{cases} \Delta g=w\cdot \na v-\ve v^{-2p_0-1}-2\mu v^{-1}|\na v|^2- \frac{1}{|\O|}\int ( w\cdot \na v-\ve v^{-2p_0-1}-2\mu v^{-1}|\na v|^2) dx,& x \in \O, \\ \int gdx=0.\end{cases} \ee
Since \eqref{zqbb89}     implies  \be \label{qbb1`}\left| \int w\cdot \na v dx \right| \le C\|w\|_{L^2 }\|\na v\|_{L^2 }\le C,\ee  we obtain   that    $\na g$  satisfies for any $p> 2$,
\be \label{qbb90} \ba\|\na g\|_{L^p }&\le C\|\Delta g\|_{L^{3p/(p+3)}}\\
&\le C(p)\|w\|_{L^p }\|\na v\|_{L^3 }+C(p)\|\na v\|_{L^p }\|\na v\|_{L^3 }+C(p)\\&\le C(p)\|w\|_{L^p } +C(p)\|\na v\|_{L^p }+C(p),\ea\ee due to \eqref{qbb88}, \eqref{zqbb89}, and \eqref{qc1}.

Setting $$\ti v(x,t)\triangleq v(x,t)+\frac{1}{2|\O|}\int_0^t\int ( w\cdot \na v-\ve v^{-2p_0-1}-2\mu v^{-1}|\na v|^2) dxdt,$$
one deduces from \eqref{zqbb87} that \be \label{qbv87} \begin{cases} 2\ti v_t-\ve\div (|\na \ti v|^2\na \ti v)=\div \ti f,\\ \ti v(x,0)=v(x,0), \end{cases}\ee with $\ti f\triangleq 2\mu  \na\ti v-wv-\na g.$

  Thus,     applying the   $L^p$-estimates  \cite[Theorem 1]{am1} (see also\cite{bog1,bor1})    to   \eqref{qbv87}  with periodic data   yields   that for any $p\ge 4 $
\be \label{qbb96}  \ba \int_0^T\|\na v \|_{L^{3p} }^{3p}dt &=\int_0^T\|\na \ti v \|_{L^{3p} }^{3p}dt \\ &\le C(p)\left(1+\int_0^T \|\ti f\|_{L^p }^p dt\right)^2\\ &\le C(p) \left(1+ \int_0^T \|w\|_{L^p }^p dt\right)^2+C(p)\left( \int_0^T \|\na \ti v\|_{L^p }^p dt\right)^2 \\ &\le C(p)+C(p) \left(  \int_0^T \|w\|_{L^p }^p dt\right)^2+\frac{1}{2}\int_0^T\|\na   v \|_{L^{3 p} }^{3 p}dt  ,\ea\ee where we have  used  \eqref{qbb90},  \eqref{qc1}, and \eqref{zqbb89}. The combination of \eqref{zqbb89} with \eqref{qbb96} gives
\be \label{hi7} \ba &  \int_0^T\int|\na v|^{15}dxdt\le C.\ea\ee

Next, it follows from \eqref{w=u0}$_2$ that
 \be \label{hi2} \ba &    w_t -(\nu+\sqrt{\ve})\Delta w- (\nu-\mu)\na \div w=F \ea\ee
with
 \be \label{hi3} \ba
F\triangleq& -   w\cdot\na w-\n^{-1}\na P +(\nu+\mu+\sqrt{\ve})\n^{-1}\na\n \cdot\na w +(\nu-\mu)\n^{-1}\na w \cdot\na \n\\
& +\ve \n^{-1}v |\na v|^2\na v\cdot\na w- \ve^{3/2}     |w|^{3}u  -r_0\n^{-1}u-r_1 |u|^2u  -\ve\n^{-p_0-1} w. \ea\ee
Multiplying \eqref{hi2} by  $|w|^2w$ and integrating the resulting equality by parts, it holds that
\be \label{hi8} \ba
&\frac{1}{4}\frac{d}{dt}\|w\|_{L^4}^4+(\nu+\sqrt{\ve})\int \xl(|\na w|^2|w|^2+\frac{1}{2}|\na|w|^2|^2\xr)dx\\
&\quad\quad\quad\quad\quad +(\nu-\mu)\int |\div w|^2|w|^2dx+\ve\int \n^{-p_0-1}|w|^4dx\\
&=-(\nu-\mu)\int \div w \na|w|^2\cdot wdx - \int w\cdot\na w\cdot w|w|^2dx-\int\n^{-1}\na P \cdot w|w|^2dx \\
&\quad+(\nu+\mu+\sqrt{\ve})\int\n^{-1}\na\n \cdot\na w\cdot w|w|^2 dx+(\nu-\mu)\int\n^{-1}|w|^2 w \cdot\na w \cdot\na \n dx\\
&\quad+\ve \int\n^{-1}v |\na v|^2\na v\cdot\na w \cdot w|w|^2dx- \ve^{3/2}    \int |w|^{3}u \cdot w|w|^2dx\\
&\quad  -r_0\int\n^{-1}u\cdot w|w|^2dx-r_1 \int|u|^2u\cdot w|w|^2dx \\
&\triangleq \sum_i^{9} J_i. \ea\ee

The straight arguments together with \eqref{qc1}, \eqref{zqbb89}, and \eqref{hi7} derive the   estimates on each $J_i(i=1,2,\cdots,10)$ as follows:
\be \label{hi9} \ba
J_1
\le \frac{3(\nu-\mu)}{4}\int |\div w|^2|w|^2dx+\frac{\nu-\mu}{3}\int |\na|w|^2|^2dx,\ea\ee
\be \label{hi10} \ba
J_2+J_3+J_8 &\le   \int |\na w| |w|^4dx +C\int |\na v| |w|^3dx+C\int |u| |w|^3dx\\
&\le  \delta\int  |\na w|^2|w|^2dx+C\int |w|^6dx+C\|\na v\|_{L^2}^2+C\|u\|_{L^2}^2\\
&\le  \delta\int  |\na w|^2|w|^2dx+C\|w\|_{L^5}^5+\frac{\ve^{3/2}}{8}\|w\|_{L^7}^7+C, \ea\ee
\be \label{hi12} \ba
J_4+J_5+J_6&\le C \int |\na w||\na v| |w|^3dx+C \int |\na w||\na v|^3 |w|^3dx\\
&\le  \delta\int |\na w|^2 |w|^2dx+C \int |\na v|^2|w|^{4}dx+C \int |\na v|^6|w|^{4}dx\\
& \le \delta\int |\na w|^2 |w|^2dx+\frac{\ve^{3/2}}{8}\|w\|_{L^7}^7+ C\|\na v\|_{L^{14}}^{14}+C\|\na v\|_{L^{14/3}}^{14/3}\\
& \le \delta\int |\na w|^2 |w|^2dx+\frac{\ve^{3/2}}{8}\|w\|_{L^7}^7+  C\|\na v\|_{L^{15}}^{15}+C,\ea\ee

\be \label{hi20} \ba
J_{7} &=-\ve^{3/2} \int  |w|^5(w-\mu\na\log\n)\cdot wdx\\
&=-\ve^{3/2} \int  |w|^7dx +\ve^{3/2}\mu \int  |w|^5 \na\log\n\cdot wdx\\
&\le-\ve^{3/2} \int  |w|^7dx +\frac{\ve^{3/2}}{2}  \int  |w|^7dx+  C\|\na v\|_{L^{15}}^{15}+C,\ea\ee
and
\be \label{hi18} \ba
J_{9}&= -\frac{r_1}{2}\int|u|^2(w-\mu\na\log\n)\cdot w|w|^2dx -\frac{r_1}{2}\int|u|^2u\cdot (u+\mu\na\log\n)|w|^2dx \\
&= -\frac{r_1}{2}\int|u|^2|w|^4dx -\frac{r_1}{2}\int|u|^4|w|^2dx  +\frac{r_1}{2}\mu^2\int|u|^2 |\na\log\n |^2|w|^2dx\\
&\le  -\frac{r_1}{2}\int|u|^2|w|^4dx -\frac{r_1}{2}\int|u|^4|w|^2dx  +\frac{r_1}{4} \int|u|^4|w|^2dx+C \int|\na v|^4|w|^2dx\\
&\le  -\frac{r_1}{2}\int|u|^2|w|^4dx -\frac{r_1}{2}\int|u|^4|w|^2dx +\frac{r_1}{4} \int|u|^4|w|^2dx+C\|w\|_{L^5}^5+ C\|\na v\|_{L^{15}}^{15}+C.
 \ea\ee

 Substituting  \eqref{hi9}--\eqref{hi18} into \eqref{hi8} and choosing $\delta$ suitably small enough, we get
\be \label{hi21} \ba
&\frac{1}{4}\frac{d}{dt}\|w\|_{L^4}^4+\frac{ \nu+\sqrt{\ve} }{4}\int |\na w|^2|w|^2 dx+\frac{\nu-\mu}{6}\int |\div w|^2|w|^2dx\\
&\quad+\frac{\ve^{3/2}}{4}  \int  |w|^7dx+\ve \int \n^{-p_0-1} |w|^4dx +\frac{r_1}{2}\int|u|^2|w|^4dx +\frac{r_1}{4}\int|u|^4|w|^2dx\\
&\le C\|w\|_{L^5}^5+ C\|\na v\|_{L^{15}}^{15}+C, \ea\ee
which together with \eqref{zqbb89}, \eqref{hi7}, and \eqref{qbb96} gives that
\be \label{hi22} \ba
\sup_{0\le t\le T}\|w\|_{L^4}^4&+\int_0^T\int \xl(|\na w|^2|w|^2 +|\div w|^2|w|^2 +  |w|^7+|\na v|^{21}\xr)dxdt\\
& +\int_0^T\int\xl(\n^{-p_0-1} |w|^4+|u|^2|w|^4 + |u|^4|w|^2\xr) dxdt\le C.  \ea\ee

Hence, \eqref{hi0} is deduced directly from \eqref{zqbb89} and \eqref{hi22}. The proof of Lemma \ref{hiu} is finished. \hfill$\Box$

  In order to obtain the   global strong solutions of   problem  \eqref{zqba1}--\eqref{zqbb35}, we still need to derive some necessary higher order   estimates on $(\n,w)$ in the  following lemma.

\begin{lemma}\label{zqlem13'}
For any $p> 2,$ there exists some constant $C $ depending on $\ve$, $r_0$,  $r_1$, $\kappa$, and $p $  such that \be\label{zqbc4}   \int_0^T\left(\|(\n_t,\na \n_t,u_t,w_t\|^p_{L^p }+\|(\n,\na \n, u,w)\|_{W^{2,p} }^p\right)dt\le C . \ee
\end{lemma}
{\it Proof.}  Multiplying  \eqref{hi2} by $- 2\Delta w $ and integrating the resulting equality over $\O$ lead to\be\ba \label{lhi1}  &\frac{d}{dt} \|\na w\|_{L^2}^2 +\int 2\xl((\nu+\sqrt{\ve})| \Delta w|^2+ (\nu-\mu)|\na\div w|^2\xr)dx \\&=-2\int\xl(-w\cdot\na w-\n^{-1}\na P +( \mu+\nu+\sqrt{\ve})\n^{-1}\na\n \cdot\na w +(\nu-\mu)\n^{-1}\na w\cdot\na\n \right.\\
&\qquad\left.+\ve \n^{-1}v |\na v|^2\na v\cdot\na w- \ve^{3/2} |w|^{3}u  -r_0\n^{-1}u-r_1 |u|^2u  -\ve\n^{-p_0-1} w\xr)\cdot \Delta wdx\\
&\triangleq \sum_{i=1}^{9} \tilde{J}_i.\ea\ee

 Using \eqref{qc1} and \eqref{hi0}, the terms $\tilde{J}_i(i=1,2,\cdots,9)$ in \eqref{lhi1} can be estimated as follows:
  \be \label{lhi2} \ba
\tilde{J}_1+\tilde{J}_3+\tilde{J}_4+\tilde{J}_5
&\le C\|\Delta w\|_{L^2}\xl(\|w\|_{L^5}+\|\na v\|_{L^5}+\|\na v\|_{L^{15}}^3\xr)\|\na w\|_{L^{10/3}}\\
&\le C\|\Delta w\|_{L^2}\xl(\|w\|_{L^5}+\|\na v\|_{L^5}+\|\na v\|_{L^{15}}^3\xr)\|\na w\|_{L^2}^{2/5}\|\na^2w\|_{L^2}^{3/5}\\
&\le \delta\|\Delta w\|_{L^2}^2+C(\|w\|_{L^5}^5+\|\na v\|_{L^5}^5+\|\na v\|_{L^{15}}^{15})\|\na w\|_{L^2}^2,\ea\ee
 \be \label{lhi3} \ba\tilde{J}_2+\sum_{i=7}^{9}\tilde{J}_i
 &\le C\int \xl( |\na v|+|u|+|u|^3+|w|\xr)| \Delta w|dx\\
 &\le \delta\|\Delta w\|_{L^2}^2 +C \|\na v\|_{L^2}^2+C \|u\|_{L^2}^2+C \|w\|_{L^2}^2+C \|u\|_{L^6}^6\\
 &\le \delta\|\Delta w\|_{L^2}^2+C+C \|w\|_{L^6}^6+C \|\na v\|_{L^6}^6\\
 &\le \delta\|\Delta w\|_{L^2}^2+C+C \|w\|_{L^7}^7+C \|\na v\|_{L^{15}}^{15}, \ea\ee
and
 \be \label{lhi6} \ba
\tilde{J}_6& =2\ve^{3/2}\int |w|^{3}(w-\mu\na\log\n)\cdot  \Delta wdx\\
&=-2\ve^{3/2}\int |w|^{3}|\na w|^2dx-6\ve^{3/2}\int |\na|w||^{2}|w|^3dx-2\ve^{3/2}\mu\int |w|^{3} \na\log\n \cdot  \Delta wdx\\
&\le -2\ve^{3/2}\int |w|^{3}|\na w|^2dx-6\ve^{3/2}\int |\na|w||^{2}|w|^3dx+C\|w\|_{L^7}^{7}+C \|\na v\|_{L^{14}}^{14}+\delta\|\Delta w\|_{L^2}^2\\
&\le -2\ve^{3/2}\int |w|^{3}|\na w|^2dx-6\ve^{3/2}\int |\na|w||^{2}|w|^3dx+C\|w\|_{L^7}^{7}+C \|\na v\|_{L^{15}}^{15}+C+\delta\|\Delta w\|_{L^2}^2. \ea\ee

Submitting  \eqref{lhi2}--\eqref{lhi6} into \eqref{lhi1}, one gets after  choosing $\delta$ suitably small enough  that
\be\ba \label{lhi11}  &\xl(\|\na w\|_{L^2}^2 \xr)_t+\int \xl((\nu+\sqrt{\ve})| \Delta w|^2+(\nu-\mu)|\na\div w|^2\xr)dx \\
&\quad+2\ve^{3/2}\int |w|^{3}|\na w|^2dx+6\ve^{3/2}\int |\na|w||^{2}|w|^3dx \\
&\le C\xl(\|w\|_{L^5}^5+\|\na v\|_{L^5}^5+\|\na v\|_{L^{15}}^{15}\xr)\|\na w\|_{L^2}^2+C\|w\|_{L^7}^{7}+C \|\na v\|_{L^{15}}^{15}+C,\ea\ee
which together with \eqref{hi0} and Gronwall's inequality yields
\be \label{qq11}\sup_{0\le t\le T}\|\na w\|^2_{L^2}+\int_0^T \|\na^2w\|_{L^2}^2dt\le C.\ee

It thus follows from \eqref{qq11} and  Sobolev inequality that \bnn \|w\|_{L^{10}(\O\times (0,T))}+\|\na w\|_{L^{10/3}(\O\times (0,T))}\le C.\enn This along with   \eqref{qbb96}--\eqref{hi3} and \eqref{qq11} gives \be \label{qq12}\|w_t\|_{L^{2}(\O\times(0,T))}+ \|\na^2 w \|_{L^{2 }(\O\times(0,T))} +\|F\|_{L^{5/2}(\O\times (0,T))} \le C.\ee
Using \eqref{qq12} and applying the standard $L^p$-estimates  to \eqref{hi2} \eqref{hi3}  \eqref{zqbb35}   with periodic data yield     that for any $p\ge 2$,
  \be  \label{qbb99} \ba  \|w_t\|_{L^{p}(\O\times(0,T))}+ \|\na^2 w \|_{L^{p }(\O\times(0,T))} \le C(p)+C(p)\|F\|_{L^{p}(\O\times(0,T))}  .  \ea\ee
In particular, the combination of    \eqref{qq12} with \eqref{qbb99}  shows
  \be\notag  \label{qbb99'} \ba  \|w_t\|_{L^{5/2}(\O\times(0,T))}+ \|\na^2 w \|_{L^{5/2 }(\O\times(0,T))} \le C  .  \ea\ee
This combined with \eqref{hi0} and the Sobolev inequality (\cite[Chapter II (3.15)]{la1}) yields that for any $q>2$,
 \bnn \ba \|  w \|_{L^{q}(\O\times(0,T))}
 +\|\na w \|_{L^{5}(\O\times(0,T))}\le C(q), \ea\enn
which along with \eqref{qbb96} and \eqref{hi3}  gives
\bnn \|F\|_{L^{9/2}(\O\times (0,T))} \le C.\enn
Combining this with \eqref{qbb99}  leads to
  \bnn   \ba  \|w_t\|_{L^{9/2}(\O\times(0,T))}+ \|\na^2 w \|_{L^{9/2 }(\O\times(0,T))} \le C ,  \ea\enn
which together with the Sobolev inequality (\cite[Chapter II (3.15)]{la1}) shows
    \bnn \ba \|w \|_{L^{\infty}(\O\times(0,T))}+\|\na w \|_{L^{45}(\O\times(0,T))}\le C . \ea\enn
Thus, we get    \be \notag\label{qq13}\|F\|_{L^{40}(\O\times (0,T))} \le C,\ee which along with \eqref{qbb99} gives
  \bnn    \ba  \|w_t\|_{L^{40}(\O\times(0,T))}+ \|\na^2 w \|_{L^{40 }(\O\times(0,T))} \le C . \ea\enn  The  Sobolev inequality (\cite[Chapter II (3.15)]{la1}) thus implies \bnn\|\na w\|_{L^{\infty}(\O\times(0,T))}\le C.\enn Then, it holds that for any $p>2 , $
 \be \label{qc3} \|w_t\|_{L^p(\O\times(0,T))}+ \|\na^2 w \|_{L^p(\O\times(0,T))} \le C(p).\ee
With \eqref{qc3}   in hand, one can deduce easily from \eqref{hi6}  and \eqref{zqbb35} that for any $p> 2,$
 \be \label{qc3-1}\|\n_t \|_{L^p( 0,T,W^{1,p}(\O))}+ \|\na^2 \n \|_{L^p( 0,T,W^{1,p}(\O))} \le C(p).\ee
Recalling the definition of $w$ in \eqref{lw}, the combination of \eqref{qc3} with \eqref{qc3-1} yields
  \be \ba\label{qc3-2} &\|u_t\|_{L^p(\O\times(0,T))}+ \|\na^2 u \|_{L^p(\O\times(0,T))}\\
  &\le \|w_t\|_{L^p(\O\times(0,T))}+ \|\na^2 w \|_{L^p(\O\times(0,T))}+C \|\na\n_t\|_{L^p(\O\times(0,T))}+ C\|\na^3 \n \|_{L^p(\O\times(0,T))}\\
  &\le C(p),\ea\ee
 which  together with \eqref{qc3}--\eqref{qc3-1}  gives the desired estimate \eqref{zqbc4},  and thus finishes the proof of Lemma \ref{zqlem13'}. \hfill$\Box$

 \section{Compactness results}

Let \be\label{qlk1}\si_0\triangleq 10^{-10},\ee
we choose
\be\nonumber 0
\le\ti \n_{0\ve}\in C^\infty (\O),\quad
\|\na\sqrt{\ti\n_{0\ve}} \|_{L^4}^4\le \ve^{-4\si_0}
\ee
satisfying
\be\nonumber\label{qpd17}
\| r_0\log_{-}\ti\n_{0\ve}-r_0\log_{-}\n_0\|_{L^1}+\| \ti\n_{0\ve}-\n_0\|_{L^1}+ \| \ti\n_{0\ve}-\n_0\|_{L^\ga}
+\|\na( \sqrt{\ti\n_{0\ve}}-\sqrt{\n_{0}})\|_{L^2}<\ve.
\ee
Set
\bnn\label{qpd7}
\n_{0\ve}=\left(\ti\n_{0\ve}^6+\ve^{24\si_0  }\right)^{\frac16},
\enn
it is easy to check that
\be\label{qpd8} \lim_{\ve\rightarrow 0}\|\n_{0\ve}  -\n_0\|_{L^1 }=0\ee and that there exists some constant $C$ independent of $\ve$ such that \eqref{qpd9} holds. Furthermore, we choose $\ti m_{0\ve}$ such that
\begin{equation*}
\|\ti m_{0\ve}-\n_0^{-1/2}m_0 \|_{L^2 }\leq\varepsilon.
\end{equation*}
Then, define $u_{0\ve}$ as follows,
\be \la{pd21} u_{0\ve}=  \n_{0\ve}^{-1/2}\ti m_{0\ve}, \ee
we thus have
\be\la{pd13}\lim_{\ve\rightarrow 0}\|\n_{0\ve} u_{0\ve} -m_0\|_{L^1 }=0.\ee  Moreover,
it is easy to check that  \eqref{qpd22}  is still valid for $(\n_{0\ve},u_{0\ve} )$.

Extending   $(\n_{0\ve},u_{0\ve} )$ $\O$-periodically to $\r^3$,
we will consider  the problem \eqref{w=u0} with the initial data  $(\n_{0\ve},w_{0\ve})$ for $w_{0\ve}\triangleq u_{0\ve}+\mu\na\log\n_{0\ve}$. The standard parabolic theory \cite{la1}  together with Lemmas \ref{zqlem13} and \ref{hiu}--\ref{zqlem13'} illustrates that there is a unique strong solution
$(\n_\ve,w_\ve)\in C([0,T),W^{2,p}(\O))$ for any $T>0$ and any $p>2$. Then,  this in turn implies that the problem \eqref{zqba1}--\eqref{zqbb35} has a unique strong solution
$(\n_\ve,u_\ve)$ such that for any $T>0$ and  any $p>2$,
\bnn  \nv,\,\uv,\,(\nv)_t ,\,\na(\nv)_t , \,(\uv)_t ,\,\na^2 \nv  ,\,\na^3 \nv  ,\,\na^2 \uv  \in L^p(\b\times(0,T)). \enn
Moreover,  all estimates obtained in  Lemmas \ref{zqlem10} and \ref{zqlem11} still hold for the solution  $(\nv,\uv) $ to the problem \eqref{zqba1}--\eqref{zqbb35}.

Letting $\ve\to 0^+,$ we will prove that $(\nv,\sqrt{\nv}\uv) $ converges, up to the extraction of subsequences, to the limit $(\n ,\sqrt{\n }u)$ in some sense. These convergences, see Lemmas \ref{qlema2}--\ref{lem36}, are crucial to show that $(\n ,\sqrt{\n }u)$ is a weak solution to \eqref{qns}--\eqref{1.2}. The proof of Lemmas \ref{qlema2}--\ref{lem36} are similar as those in Li-Xin \cite{li04} (see also partially in \cite{vy1,AS2017}), which are sketched here for completeness.

 We begin with the following strong convergence of $\sqrt{\nv} $ and $\nv.$

\begin{lemma} \label{qlema2}   There exists a function $\n\in L^\infty(0,T; L^1\cap L^\ga )$ such that up to a subsequence,
 \be \label{qee3}\sqrt{\nv} \ro \sqrt{\n} \ \ \mbox{strongly in} \ \ L^2(0,T;H^1 ),\ee
 \be \label{qe3}\nv \ro \n \ \ \mbox{strongly   in }\ \ L^\ga (\O\times(0,T )),\ee
  \be \label{gyqe3}\na^2\sqrt{\nv} \ro \na^2\sqrt{\n} \ \ \mbox{weakly   in }\ \ L^2  (\O\times(0,T )).\ee
 In particular, it holds
 \be \la{ze3}\nv  \ro \n  \mbox{ almost everywhere  in  } \O\times(0,T  ).\ee   \end{lemma}
  {\it Proof.} First, for   $\vf\triangleq \sqrt{\nv},$ it follows from \eqref{zqbb30}, \eqref{zqbb5'}, \eqref{useful}, and \eqref{newb}     that there exists some generic positive constant $C$ independent of $\ve  $      such that
\begin{align}\label{qd1}
& \sup_{0\le t\le T}\int  (\n_\ve |u_\ve|^2+\n_\ve+\n_\ve^\ga+\ve\nv^{-p_0})dx    +\ioo|\na(\sqrt{\n_\ve}u_\ve)-u_\ve\otimes \na\sqrt{\n_\ve}|^2dxdt\nonumber \\
& +\ve\ioo \left( |\na\vf|^4|\uv|^2+\nv^{-p_0}|\uv|^2
+\ve^{1/2}\nv|w_\varepsilon|^3|u_\varepsilon|^2
+\ve^{1/2}\nv|w_\varepsilon|^5+\ve^{1/2}\nv |u_\varepsilon|^5\right) dxdt \nonumber \\
&+\ioo \xl(\n_\ve |\na u_\ve|^2 +\n_\ve|\uv|^{4} + |\uv|^{2}\xr)dx dt+ \ve^{2} \ioo \nv^{-2p_0-1}dxdt
\le C,
\end{align}
and
\begin{align}\label{qd2}
    & \sup_{0\le t\le T} \int (|\na\vf|^2+\ve|\na\vf|^4) dx   + \ioo  \left( |\na \n_\ve^{\ga/2}|^2+\n_\ve |\nabla^2\log\n_\ve|^2+|\nabla^2 v_\ve|^2+|\na \nv^{\frac{1}{4}}|^4\right.\nonumber \\
    &\left.+|\na (\sqrt{\nv}\uv)|^2
\right)dxdt  +\ve \ioo   \left(v_\ve^{-2}|\na v_\ve|^6+v_\ve^{-3}|\na v_\ve|^5+|\na\vf|^2|\na^2\vf|^2 +\ve|\na\vf|^4|\na^2\vf|^2 \right)  dxdt  \le C.
\end{align}
Then, one deduces  from \eqref{qd1}, \eqref{qd2}, H{\"o}lder  and Sobolev  inequalities that
\be \ba \label{qq7}
\ve^{\frac43}\int_0^T\|\na \vf\|_{L^6}^6dt
&\le \int_0^T\|\na \vf\|_{L^2}^{\frac23}(\ve\|\na \vf\|^4_{L^4})^{\frac13}(\ve\|\na \vf\|_{L^{12}}^4)dt   \\
&\le C\ve\int_0^T\||\na \vf||\na^2\vf|\|_{L^2}^2dt   \le C.
\ea\ee

Since $\nv $ satisfies
\be \label{qe5}\ba (\n_\ve )_t+\div (\nv u_\ve)  = {\ve} \vf\div(|\na\vf|^2\na \vf)+\ve\nv^{-p_0}, \ea\ee
by assuming $\nv>0,$ we may rewrite \eqref{qe5} as follows
\be \label{qee5}\ba 2(\sqrt{\n_\ve})_t  = -2\div (\sqrt{\nv} u_\ve)
+\sqrt{\nv} \div u_\ve
+{\ve} \div(|\na\vf|^2\na \vf)+\ve\nv^{-p_0-\frac12}.
\ea\ee
It follows from \eqref{qd1}, \eqref{qd2}, and \eqref{qq7} that
\be \label{qe7}\ba \ve^2\ioo |\na\vf|^6dxdt \le C \ve^{\frac23} \ea\ee
and \be \label{qe7'}\ba   \ioo \left(\n_\ve  |u_\ve|^2+\n_\ve(\div u_\ve)^2\right)dxdt+\ve^2\ioo \nv^{-2p_0-1}dxdt \le C.\ea\ee
The combination of \eqref{qee5}--\eqref{qe7'} implies that \be\label{qve4}\|(\sqrt{\nv})_t\|_{L^2(0,T;H^{-1} )}\le C.\ee

Furthermore, it is easy to derive  from   \eqref{qd1} and \eqref{qd2} that
\be\label{qvee4}
\|\sqrt{\nv}\|_{L^2(0,T;H^{2} )}\le C,
\ee
which combined with \eqref{qve4} and Aubin-Lions lemma yields \eqref{qee3}.

Next, we claim that for $\gamma\in (1,3)$,
\begin{equation}\label{3.019}
\|\n_\ve^{\gamma}\|_{L^{\frac53}((0,T)\times\Omega)}\leq C.
\end{equation}
This along with \eqref{qee3} yields directly the desired \eqref{qe3} and \eqref{ze3}. Furthermore, the convergence  \eqref{gyqe3} is deduced  directly form \eqref{qvee4} and \eqref{qe3}.

Now, it remains  to prove \eqref{3.019}.  It is easy to deduce  from \eqref{qd2} that
\begin{equation*}
\|\nabla\n_\ve^{\frac{\gamma}{2}}\|_{L^2((0,T)\times\Omega)}\leq C,
\end{equation*}
which together with Sobolev's embedding theorem gives
\begin{equation}\label{3.19}
\|\n_\ve^{\gamma}\|_{L^1(0,T;L^3 )}\leq C.
\end{equation}
Note that \eqref{qd1} implies that
$$\|\n_\ve^\gamma\|_{L^\infty(0,T;L^1 )}\leq C,$$
this combined with \eqref{3.19} yields  \eqref{3.019}.  The proof of Lemma \ref{qlema2} is finished.
\hfill $\Box$

Next, we have the following lemma which  deals with the compactness of the momentum.

\begin{lemma}\label{qvlem03} There exists a function $m(x,t)\in L^2(0,T;L^{\frac32} )$ such that up to a subsequence,
\be \label{qve10}
\nv\uv\rt m \mbox{ in }L^2(0,T;L^p )
\ee
for all $p\in [1,\frac32)$.
Moreover, there exists a function $u$ in $L^2((0,T)\times \Omega)$ such that up to a subsequence
\begin{align}\label{velocity}
\uv \to u \text{ weakly in } L^2((0,T)\times \Omega).
\end{align}
And, it holds that\be \label{qve12}\nv\uv\rt \n u \mbox{ almost everywhere }(x,t)\in\O\times (0,T) .\ee
 \end{lemma}
{\it Proof.}
First, it follows from H{\"o}lder inequality, \eqref{qd1}, and \eqref{qd2} that
\begin{align}\label{qve11}
\int_0^T\|\na (\nv \uv)\|_{L^1 }^2dt
& \le C\int_0^T\left( \|\nv \|_{L^{1} }\|\sqrt{\nv}\na \uv\|_{L^2 }^2
+ \|\sqrt{\nv} \uv\|_{L^2}^2 \|\na\sqrt{\nv}\|_{L^2 }^2\right) dt \le  C
\end{align}
and
\begin{align*}
\sup_{0\le t\le T} \|\nv \uv\|_{L^1} \le C\sup_{0\le t\le T} (\|\nv\|_{L^1} + \|\nv |\uv|^2 \|_{L^1}) \le C.
\end{align*}
Hence, one has
\begin{align*}
\|\nv \uv \|_{L^2(0,T; W^{1,1})}\le C.
\end{align*}

Next, the straight calculations show that
\be\ba\label{qve14}
 &(\nv\uv)_t+\div(\nv\uv\otimes\uv)-2\nu\div(\nv \mathcal{D}\uv)
 -\ka^2\div\left(\nv\nabla^2\log\nv\right)+\na P(\nv)  \\
 &=\ve\div( \vf|\na\vf|^2\na\vf\otimes \uv)-\ve |\na\vf|^4\uv -\ve\mu \na\nv^{-p_0} +\sqrt{\ve}\div( \nv  \na\uv )  \\
 &\quad -\ve^{\frac32}\nv|w_\varepsilon|^3\uv-r_0\uv-r_1\nv|\uv|^2\uv+\sqrt{\ve}\mu\div(\nv\na^2\log\nv)
-\ve\mu |\na v_\varepsilon|^4\na \log\nv
  \\
 &\quad -\ve\mu\na( v_\varepsilon\div( |\na v_\varepsilon|^2 \na v_\varepsilon))
+\ve\mu\div( v_\varepsilon|\na v_\varepsilon|^2 \na v_\varepsilon\otimes\na \log \nv).
\ea\ee

For the terms on the left-hand side of \eqref{qve14}, one has
\be\label{qve15}\ba
& \ioo  \nv|\uv|^2dxdt  +\ioo  \nv^\ga dxdt\le C,
\ea\ee
\be \label{qve18}\ba
\ioo\nv|\na\uv|dxdt
& \le C  \ioo \nv |\na\uv|^2 dxdt+C\ioo  \nv  dxdt   \le C,
\ea\ee
\begin{align}\label{zqve18}
\ioo\nv |\na^2\log\nv|dxdt
\le C\ioo \nv dxdt+C\ioo\nv|\na^2\log\nv|^2dxdt
\le C.
\end{align}

Using \eqref{qd1} and \eqref{qd2}, we can estimate each term on the right-hand side of \eqref{qve14} as follows:
\be \label{qqve1}\ba
& \ve\ioo\left(\vf |\na\vf |^3| \uv|+ |\na\vf|^4|\uv|\right)dxdt \\
& \le C \ve  \int_0^T\|\vf\uv\|_{L^2}\|\na \vf\|^3_{L^6}dt   +C\left(\ve\ioo  |\na\vf|^4|\uv|^2dxdt \right)^{\frac12}\left(\ve\ioo  |\na\vf|^4 dxdt \right)^{\frac12}   \\
& \le C\ve\xl(\ioo |\na\vf|^6dxdt\xr)^{1/2}+C\left[\ve\xl(\ioo |\na\vf|^6dxdt\xr)^{1/2} \right]^{\frac12}\\
& \le C\ve^{\frac16},
\ea\ee
where in the last inequality one has used \eqref{qq7}.
Moreover, it holds
\begin{align} \la{eqe7'}
\ve\ioo \nv^{-p_0}dxdt & \le \ve^{\frac{1}{2p_0+1}}\left(\ve^{2}\ioo \nv^{-2p_0-1}dxdt\right)^{\frac{p_0}{2p_0+1}}  \le C \ve^{\frac{1}{2p_0+1}},
\end{align}
and
\be \label{zqqve1}\ba
& \ioo(|\uv|+\nv|\uv|^3)dxdt   \\
& \le C\xl(\int_0^T\int|\uv|^2dxdt\xr)^{\frac12}+C\left(\ioo \nv|\uv|^2dxdt \right)^{\frac12}\left(\ioo \nv|\uv|^4 dxdt \right)^{\frac12}   \\
& \le C.
\ea\ee
The  H{\"o}lder  inequality together with \eqref{qd1}  and \eqref{qd2} yields
\be\ba \label{zqqve2}
&\ve^{\frac32}\ioo\nv|w_\varepsilon|^3|\uv|dxdt\\
& \le C \ve^{\frac32}\xl(\int_0^T\int\nv|w_\varepsilon|^5dxdt\xr)^{3/5}\xl(\int_0^T\int\nv|\uv|^{5/2}dxdt\xr)^{2/5}  \\
& \le C\ve^{\frac32}\xl(\int_0^T\int\nv|w_\varepsilon|^5dxdt\xr)^{3/5}\xl(\xl(\int_0^T\int\nv|\uv|^5dxdt\xr)^{1/2}\xl(\int_0^T\int\nv dxdt\xr)^{1/2}\xr)^{2/5}  \\
& \le C\ve^{\frac{3}{10}}.
\ea\ee
It follows from \eqref{qd1}, \eqref{qd2}, \eqref{qq7}, H{\"o}lder  and  Sobolev inequalities that
\begin{align}\label{qqev1}
\ve\ioo\vf|\div(|\na\vf|^2\na\vf)|dxdt
& \le C\ve\int_0^T\int \vf|\na\vf|^2|\na^2\vf| dxdt \nonumber \\
& \le C\ve\int_0^T\|\vf\|_{L^6}\|\na\vf\|_{L^6}^2\|\na^2\vf\|_{L^2}dt \nonumber \\
& \le C\sup_{0\le t\le T}(\|\na \vf\|_{L^2})\ve\left(\int_0^T\|\na\vf\|_{L^6}^4dt\right)^{\frac12}
\left(\int_0^T\|\na^2\vf\|_{L^2}^2dt\right)^{\frac12}
\nonumber \\
& \le C\ve^{\frac59}\left(\ve^{\frac43}\int_0^T\|\na\vf\|_{L^6}^6dt\right)^{\frac13}
\nonumber \\
& \le C\ve^{\frac59}.
\end{align}
Finally, we deduce from H{\"o}lder  inequality and  \eqref{qd2}  that
\be\ba\label{qqev3}
&\ve\ioo|\na\vf|^4|\nabla\log\nv|dxdt+\ve\ioo\vf|\na\vf|^3|\nabla\log\nv|dxdt\\
&\le C\ve\ioo\vf^{-1}|\na\vf|^5dxdt+  C\ve\ioo|\na\vf|^4dxdt   \\
&\le C\ve^{7/6} \ioo \vf^{-2}|\na\vf|^6dxdt+  C\ve^{5/6}\ioo|\na\vf|^4dxdt   \\
& \le C\ve^{1/6}+C\ve^{5/6}\int_0^T\|\na\vf\|_{L^2}\|\na\vf\|_{L^6}^3dt\\
&\le C\ve^{1/6}.
\ea\ee
The combination of \eqref{qd1}--\eqref{qd2} with \eqref{qve14}--\eqref{qqev3} leads to
\be \label{qnk5}
\|(\nv\uv)_t\|_{ L^1(0,T;W^{-1,1} )}\le C.
\ee
Hence, \eqref{qve10} is deduced from  Aubin-Lions lemma,  \eqref{qve11}, and \eqref{qnk5}.

Next, it's noted that $u_\varepsilon$ is uniformly bounded in $L^2((0,T)\times\O)$, which yields directly \eqref{velocity}.

Now,  it follows from \eqref{qve10}  that
\begin{align}\label{z}
\nv\uv\rt m \mbox{ almost everywhere }(x,t)\in\O\times (0,T).
\end{align}
On the one hand,  \eqref{z}   and \eqref{ze3}  show that
\begin{align}\label{x1}
\uv = \frac{\nv \uv}{\nv}\rt \frac{m}{\rho} \mbox{ almost everywhere }~\{(x,t)\in \O\times(0,T)|\n(x,t)>0\},
\end{align}
which together with  \eqref{velocity}  gives that for $\rho>0$,
\begin{align*}
m= \rho u .
\end{align*}
On the other hand, it follows from Fatou's lemma and \eqref{qd1} that
\begin{align*}
\ioo \liminf_{\ve\rt 0^+} \frac{|\rho_\ve u_\ve|^2}{\rho_\ve}dxdt =\ioo \liminf_{\ve\rt 0^+} \rho_\ve |u_\ve|^2dxdt \le \liminf_{\ve\rt 0^+} \ioo  \nv|\uv|^2dxdt\le C.
\end{align*}
This implies that if $\rho=0$, it has
\begin{align*}
m=0.
\end{align*}
Then, \eqref{qve12} is proved.
The proof of Lemma \ref{qvlem03} is completed.\hfill $\Box$

With Lemmas \ref{qlema2} and \ref{qvlem03} in hand, we are now in a position to prove the strong convergence of $\sqrt{\nv}\uv$. This is crucial for deriving the global existence of the weak solution.
\begin{lemma}\label{qlema1} Up to a subsequence, it holds
\be \label{qe42} \sqrt{\nv}\uv\rt \sqrt{\n}u \mbox{ strongly in } L^2 ( 0,T;L^2 ),\ee with \be \label{qe'42} \sqrt{\n}u \in L^\infty(0,T;L^2).\ee Moreover, it holds that
\begin{align}\label{gy061}
\sqrt{\nv}\uv\rt \sqrt{\n}u \mbox{ almost everywhere }(x,t)\in\O\times (0,T).
\end{align}\end{lemma}
{\it Proof.}
For any $M>0$, the straight calculation shows that
\be\ba \la{e47}
& \ioo|\sqrt{\nv}\uv-\sqrt{\n}u|^2dxdt  \\
& \le 2\ioo   |\sqrt{\nv}\uv1_{(|\uv|\le M)}-\sqrt{\n}u1_{(|u|\le M)}|^2dxdt \\
& \quad +2\ioo|\sqrt{\nv}\uv1_{(|\uv|\ge M)}|^2dxdt
+2\ioo|\sqrt{\n}u1_{(|u|\ge M)}|^2dxdt  \\
& \le 2\ioo   |\sqrt{\nv}\uv1_{(|\uv|\le M)}-\sqrt{\n}u1_{(|u|\le M)}|^2dxdt +\frac{2}{M^2}\ioo   \left( \nv |\uv|^4+ \n |u|^4\right)dxdt.
\ea\ee

First, it follows  from \eqref{qve12}  and  \eqref{ze3}  that
 \be\ba\label{gyg1}  \sqrt{\nv}\uv \rt \sqrt{\n}u \mbox{ almost everywhere in } ~\{(x,t)\in \O\times(0,T)|\n(x,t)>0\}. \ea\ee
 Moreover, since
\be \la{e3a}
\sqrt{\nv}|\uv|1_{(|\uv|\le M)}\le M\sqrt{\nv}
\ee
and
\be\ba \label{gyg2}  \nv  \rt  \n  \mbox{ almost everywhere in } ~\{(x,t)\in \O\times(0,T)|\n(x,t)=0\},  \ea\ee
we have
\be \nonumber\la{e45}
\sqrt{\nv}\uv1_{(|\uv|\le M)} \rt \sqrt{\n}u1_{(| u|\le M)} \mbox{ almost everywhere in } \O\times(0,T),
\ee
which, together with \eqref{e3a} and \eqref{qe3}, implies
\be\la{e48}
\lim_{\ve\rt 0^+}\ioo \left|\sqrt{\nv}\uv1_{(|\uv|\le M)}-\sqrt{\n}u1_{(|u|\le M)}\right|^2dxdt =0.
\ee

Next, Lemma \ref{zqlem10} yields that there exists
some constant $C$ independent of $\ve$ such that
\be\la{d3}\int_0^T\int \n_\ve  |u_\ve|^4dxdt\le C,\ee  which, together with  \eqref{qve12},   \eqref{ze3}, and Fatou's lemma, gives
\be\la{e44} \ba
\ioo \n|u|^4dxdt\le \ioo \liminf_{\ve\rt 0^+} \nv|\uv|^4dxdt\le \liminf_{\ve\rt 0^+} \ioo  \nv|\uv|^4dxdt\le C.
\ea\ee

Substituting  \eqref{e48}--\eqref{e44} into  \eqref{e47} yields that up to a subsequence \be \la{e51} \limsup_{\ve\rt 0^+}\ioo   |\sqrt{\nv}\uv-\sqrt{\n}u|^2dxdt\le \frac{C}{M^2},~~\mbox{for any}~M>0.\ee
We thus obtain  \eqref{qe42}  by taking $M\rt \infty$ in \eqref{e51}.
 The combination of \eqref{qd1} with  \eqref{qe42}  gives   \eqref{qe'42}.  The proof of Lemma \ref{qlema1} is finished.
\hfill $\Box$

Similar  to the proof of Lemma \ref{qlema1}, we can establish the following convergence of the damping terms.
\begin{lemma}\label{lem34}
Up to a subsequence, it holds
\begin{equation} \label{3.42}
\nv|\uv|^2\uv\rt \n |u|^2u \mbox{ strongly in } L^1(0,T;L^1 ).
\end{equation}
\end{lemma}
{\it Proof.} The direct calculation shows that for any $M>0$,
\begin{align}\la{qe47}
& \ioo |\nv|\uv|^2\uv-\n|u|^2u|dxdt \nonumber \\
& \le \ioo |\nv|\uv|^2\uv1_{(|\uv|\le M)}-\n|u|^2u1_{(|u|\le M)}|dxdt \nonumber \\
&\quad +2\ioo\nv|\uv|^31_{(|\uv|\ge M)}dxdt
+2\ioo\n|u|^31_{(|u|\ge M)}dxdt.
\end{align}

First, it follows from \eqref{x1} and \eqref{ze3} that
\begin{align}\label{gyx1}
\nv|\uv|^2\uv  \rt \n|u|^2u~\mbox{ almost everywhere in } ~\{(x,t)\in \O\times(0,T)|\n(x,t)>0\}.
\end{align}
Moreover, since
\be \la{qe3a} \nv|\uv|^2\uv1_{(|\uv|\le M)}\le M^3\nv,\ee
which together with (\ref{gyg2}) implies that
\be \nonumber\la{e45} \nv|\uv|^2\uv1_{(|\uv|\le M)} \rt \n|u|^2u1_{(|u|\le M)} \mbox{ almost everywhere in } \O\times(0,T). \ee
Then, it holds that
\be \la{qe48}
\ioo |\nv|\uv|^2\uv1_{(|\uv|\le M)}-\n|u|^2u1_{(|u|\le M)}|dxdt \rt 0\ \ \text{as}\ \ \varepsilon\rightarrow0^+.
\ee

Next, it follows from \eqref{d3} and \eqref{e44} that
\begin{align}\la{qe50}
& \ioo \xl(\nv|\uv|^31_{(|\uv|\ge M)}+\n |u|^31_{(|u|\ge M)}\xr)dxdt  \le \frac{1}{M}\ioo \left( \nv |\uv|^4+ \n |u|^4\right)dxdt  \le \frac{C}{M}.
\end{align}

Substituting \eqref{qe48} and \eqref{qe50} into \eqref{qe47} yields that up to a subsequence
\be \la{qe51}
\limsup_{\ve\rt 0^+}\ioo|\nv|\uv|^2\uv-\n|u|^2u|dxdt\le \frac{C}{M},~~\mbox{for any}~M>0.\ee
We thus obtain \eqref{3.42} by taking $M\rt \infty$ in \eqref{qe51}.
The proof of Lemma \ref{lem34} is completed.
\hfill $\Box$

Moreover, we can show the following lemma, which shows that $\na(\sqrt{\n}u)- u \otimes \na\sqrt{\n}$ is indeed a function in $L^2((0,T)\times \Omega)$ and is the limit of   $\na(\sqrt{\nv}\uv)- \uv \otimes \na\sqrt{\nv}$ in the sense of distribution.
\begin{lemma}\label{lem36}
	Up to a subsequence, it holds that
	\be \label{qe54} \na(\sqrt{\nv}\uv)- \uv \otimes \na\sqrt{\nv} \rt \na(\sqrt{\n}u)- u \otimes \na\sqrt{\n}
	\mbox{ in } \mathcal{D'}((0,T)\times\O),\ee
	\be \label{qe55} \na^{tr}(\sqrt{\nv}\uv)- \na\sqrt{\nv}\otimes \uv \rt \na^{tr}(\sqrt{\n}u)- \na\sqrt{\n} \otimes u
	\mbox{ in } \mathcal{D'}((0,T)\times\O).\ee
	Furthermore, it holds
	\begin{align}\label{new}
	\int_0^T\int |\na(\sqrt{\n}u)- u \otimes \na\sqrt{\n} |^2dxdt\le C+Cr_0+Cr_1.
	\end{align}
	%
\end{lemma}
{\it Proof.} It is easy to deduce  from (\ref{qee3}) and (\ref{velocity}) that
\begin{align}\label{z3}
 \uv\otimes \na\sqrt{\nv} \rt u \otimes  \na\sqrt{\n} \text{ in } \mathcal{D'}((0,T)\times\O),
\end{align} which together with \eqref{qe42} gives \eqref{qe54}  and thus \eqref{qe55}. Furthermore, \eqref{new} is obtained directly from   \eqref{zqbb5'} and (\ref{qe54}). This   finishes the proof of Lemma \ref{lem36}.
\hfill $\Box$

\section{Proof  of Theorem  \ref{qvth2}}

We will follow the arguments in \cite[Section 2.3]{li04} to prove that the  limit (in some sense) $(\n ,\sqrt{\n }u)$ of $(\nv,\sqrt{\nv}\uv) $ (up to a subsequence) is a weak solution to \eqref{qns}--\eqref{1.2}.

First, it follows from \eqref{qd1}, \eqref{qd2}, and \eqref{qq7} that
\begin{align} \la{eqe7}
\ve \ioo \left(\vf|\na\vf|^3+|\na \vf|^4\right)dxdt
& \le C \ve\int_0^T  \left(\|\vf\|_{L^2}\|\na\vf\|_{L^6}^3+\|\na\vf\|_{L^2}\|\na\vf\|_{L^6}^3\right)dt \nonumber \\
& \le C \ve^{\frac13}\left(\ve^{\frac43}\int_0^T\|\na\vf\|_{L^6}^6dt\right)^{\frac12} \nonumber \\
& \le C \ve^{\frac13}.
\end{align}

Then, on the one hand, for any test function $\psi$, multiplying \eqref{qe5} by $\psi$, integrating the resulting equality over $\O\times (0,T),$ and taking $\ve\rt 0$ (up to a subsequence), one can verify easily after using \eqref{qe3}, \eqref{qe42},  \eqref{qpd8}, \eqref{eqe7'}, and \eqref{eqe7}   that $(\n,\sqrt{\n}u)$ satisfies \eqref{fin1}.

On the other hand, let $\phi$ be a test function. Multiplying \eqref{qve14} by $\phi$, integrating the resulting equality over $\O\times (0,T),$ and taking $\ve\rt 0$ (up to a subsequence), by Lemmas \ref{qlema2}, \ref{qlema1}, and \ref{lem34},  
we obtain after using \eqref{qve18}--\eqref{qqve1} and \eqref{zqqve2}--\eqref{qqev3} that $(\n,\sqrt{\n}u)$ satisfies \eqref{fin2}.

The proof of Theorem \ref{qvth2} is completed.
\hfill $\Box$

 \section{Proof  of Theorem  \ref{newth}}

The  system \eqref{qns} without damping terms  is as follows:
 \be\la{imsqns}
\begin{cases}
\n_t+\div(\n u)  = 0, \\
(\n u)_t+ \div(\n u \otimes u)-2\nu\div(\n\mathcal{D} u) +\na P -2\ka^2\n\na\left(\frac{\Delta \sqrt{\n}}{\sqrt{\n}}\right)=0.
\end{cases}
\ee
We will consider the system \eqref{imsqns} on bounded domain $\O=\mathbb{T}^3$ with periodic boundary conditions and the initial conditions \eqref{1.2}. The notion of the weak solution of problem \eqref{imsqns} \eqref{1.2} is defined by $(\n,\sqrt{\n}u)$ satisfying   \eqref{fin1} and \eqref{nfin2}.

 We will consider the  approximate system of \eqref{imsqns} by choosing $r_0=r_1=0$ in \eqref{zqba1}, that is,
 \be \label{imszqba1}
 \begin{cases}
 \n_t+\div(\n u)    = \ve v\div( |\na v|^2 \na v)+\ve\n^{-p_0},
 \\   \n u_t+ \n u\cdot \na u   -2\nu\div( \n    \mathcal{D} u) +\na P -2\ka^2\n\na\left(\frac{\Delta v}{v}\right)\\=  \sqrt{\ve}\div(\n \na u)+\sqrt{\ve}\mu\div(\n \na^2\log \n) +\ve v |\na v|^2 \na v\cdot \na u+\ve\mu v |\na v|^2 \na v\cdot \na (\na \log\n)\\ \quad -\ve \n^{-p_0} u- \ve^{\frac32}\n    |w|^{3}u-\ve\mu\na \n^{-p_0} - \ve\mu\na( v\div( |\na v|^2 \na v))+\ve\mu v\div( |\na v|^2 \na v)\na \log \n.
 \end{cases}
 \ee
 In order to obtain the global existence of weak solution to the problem \eqref{imsqns} \eqref{1.2}, the main arguments here are to ensure  the smooth approximate solutions   satisfying  the a priori bounds in [3], where the compactness of finite weak solutions is shown clearly. Indeed,   one needs to prove that  the smooth  solutions to system \eqref{imszqba1} satisfying the energy estimate, the BD entropy inequality, and the Mellet-Vasseur type estimate.

It is clear that   both the  energy estimate and the BD entropy inequality obtained in Lemmas \ref{zqlem10}--\ref{zqlem11} are independent of $r_0$ and $r_1$.
 Hence, letting  $r_0=r_1=0$ in  Lemmas \ref{zqlem10}--\ref{zqlem11}, we can get directly the energy and BD entropy estimates on the smooth  solutions to system \eqref{imszqba1} as follows:

\begin{lemma} \label{imslem1} Suppose  that $11\kappa \le \nu$,   there exists some generic constant $C$ independent of $\ve$ and $\kappa$  such that
\be\label{imszqbb30}\ba& \sup_{0\le t\le T}\int\left(\n |u|^2+\n+  \n^\ga+\ve\n^{-p_0} +  \kappa^2 |\na v|^2+\ve\mu|\na v|^4 \right)dx  \\
 &+ \nu \int_0^T\int    \n   |\mathcal{D} u|^2 dxdt + \sqrt{\ve}\int_0^T\int    \n   |\na u|^2 dxdt  \\
 &+ \kappa^2 \ve\int_0^T\int   \left(   |\na v|^2|\na^2 v|^2  +|\na|\na v|^2 |^2 + |\na v|^2v^{-2p_0-1}\right) dxdt\\ & +\ve^2 \int_0^T    \int \left(\mu|\na v|^4|\na^2 v|^2+\mu|\na v|^4|\na|\na v|  |^2
 + |\na v|^4v^{-2p_0-2}+\n^{-2p_0-1}\right)dx  dt  \le C,\ea\ee
and
\be\label{imszqbb5'}\ba& \sup_{0\le t\le T}\int  \left( |\na v|^2+\ve |\na v|^4   \right)dx   + \int_0^T \int \left( \n   |\na u|^2   +   \n^{\ga-2}   |\na \n|^2\right)dxdt\\
 &+ \ka^2 \int_0^T \int \n |\na^2\log \n|^2 dxdt+\ve\nu\int_0^T\int \xl(|\na v|^2 |\na^2v|^2+|\na v|^2 |\na |\na v||^2  +\n^{-p_0-1}|\na v|^2\xr)dxdt\\&+\ve^2 \int_0^T\int \xl(|\na v|^4|\na^2v|^2+|\na v|^4 |\na |\na v||^2+\n^{-p_0-1}|\na v|^4\xr)dxdt\\
 & + \ve^{\frac{3}{2}} \int_0^T\int \xl(\rho |w|^5+ \rho |u|^5 \xr)dxdt + \ve \int_0^T\int \xl( v^{-2 } |\nabla v|^6 + v^{-3}|\nabla v|^5\xr) dxdt \le C.\ea\ee
\end{lemma}

Now, we need only to prove the Mellet-Vasseur
type estimate. Motivated by   \cite{AS2017,AS2015}, this is obtained by considering the following equivalent  transformation system of  $(\n,w)$:
\be\label{mvmu} \begin{cases}   \n_t+\div(\n w)    =\mu \Delta \n+  \ve v\div( |\na v|^2 \na v)+\ve\n^{-p_0},
     \\
     \n w_t+  \n w\cdot\na w   +\na P  -2(\nu-\mu)\div( \n    \mathcal{D} w) - \mu\n \Delta w- \sqrt{\ve}\div(\n \na w)  \\
  \quad= 2\mu\na\n \cdot\na w +\ve v|\na v|^2\na v\cdot\na w- \ve^{3/2}\n    |w|^{3}u   -\ve\n^{-p_0} w,   \end{cases}\ee
which is deduced with the same arguments as Lemma \ref{uw0}.
\begin{lemma} \label{imslem2} Suppose that  $11\kappa \le \nu$,  there exists some generic constant $C$ independent of $\ve$   such that
\be\label{mvgj}\ba& \sup_{0\le t\le T}\int \n (e+|u|^2)\ln (e+|u|^2)dx  \le C.\ea\ee
\end{lemma}
 {\it Proof.} Notice that the definition of $w$ in \eqref{lw}, one needs only to prove
 \be\label{mvgj0}\ba& \sup_{0\le t\le T}\int \n (e+|w|^2)\ln (e+|w|^2)dx  \le C.\ea\ee
 Multiplying \eqref{mvmu}$_2$ by $H\triangleq (1+ \ln (e+|w|^2)) w$ and integrating by parts yield
 \be\label{mvgj1}\ba
 &\frac{1}{2}\frac{d}{dt} \int \n (e+|w|^2)\ln (e+|w|^2)dx +\int \n  \ln (e+|w|^2) \xl(2(\nu-\mu)|\mathcal{D} w|^2+ \sqrt{\ve}\ |\na w|^2\xr)dx\\
 &\le  \frac{1}{2}\mu\int \Delta \n   (e+|w|^2)\ln (e+|w|^2)dx+2\mu\int \na \n\cdot\na w\cdot  Hdx+ \mu\int   \n \Delta w\cdot  Hdx\\
 &\quad+\frac{1}{2} \ve \int v\div(|\na v|^2\na v)   (e+|w|^2)\ln (e+|w|^2)dx+ \ve \int v |\na v|^2 \na v\cdot\na w\cdot H dx\\ &\quad+\frac{1}{2}  \ve \int \n^{-p_0}   (e+|w|^2)\ln (e+|w|^2)dx-\ve \int \n^{-p_0}  w\cdot Hdx \\
 &\quad -\int \na P\cdot Hdx-\ve^{3/2}\int \n |w|^3u\cdot Hdx  + C \int \n |\na
 w|^2dx\\
 &\triangleq \sum_{i=1}^{9}K_i + C \int \n |\na
 w|^2dx.\ea\ee

The terms $K_i(i=1,2,\cdots,9)$ in \eqref{mvgj1} can be bounded as follows.
It is easy to deduce  that
\be\label{mvgj5}\ba
 K_4+K_5& =   - \frac{\ve}{2}\int     (e+|w|^2) \ln (e+|w|^2)|\na v|^4 dx\le 0,\ea\ee
 and
 \be\label{mvgj6}\ba
 K_6+K_7& \le C\ve  \int   \n^{-p_0} dx.\ea\ee

Furthermore, integration by parts gives
\be\label{mvgj2}\ba
 K_3& =  \mu\int   \n \Delta w\cdot  (1+ \ln (e+|w|^2)) wdx\\
 &= - \mu\int   \na \n \cdot\na w\cdot  (1+ \ln (e+|w|^2)) wdx- \mu\int   \n \na \ln (e+|w|^2) \cdot\na w\cdot  wdx\\
 &\quad- \mu\int    \n |\na w|^2  (1+ \ln (e+|w|^2)) dx\\
 &=-\frac{K_2}{2}- \frac{\mu}{2}\int   \n  (e+|w|^2)^{-1}|\na |w|^2|^2 dx - \mu\int    \n |\na w|^2  (1+ \ln (e+|w|^2)) dx,\ea\ee
 and
 \be\label{mvgj3}\ba
\frac{ K_2}{2}& =  \mu\int \na \n\cdot\na w\cdot   (1+ \ln (e+|w|^2)) wdx=  \frac{\mu}{2}\int \na \n\cdot\na (e+|w|^2)\cdot   (1+ \ln (e+|w|^2)) dx\\
  & =  \frac{\mu}{2}\int \na\n\cdot\na \xl((e+|w|^2) \ln (e+|w|^2) \xr) dx  =  -\frac{\mu}{2}\int \Delta \n (e+|w|^2) \ln (e+|w|^2)  dx=-K_1.\ea\ee
The combination of \eqref{mvgj2} with \eqref{mvgj3} gives
\be\label{mvgj4}\ba
 K_1+K_2+K_3& =   - \frac{\mu}{2}\int   \n  (e+|w|^2)^{-1}|\na |w|^2|^2 dx - \mu\int    \n |\na w|^2  (1+ \ln (e+|w|^2)) dx\le 0.\ea\ee

 For the term $K_8$, it holds
 \be\label{mvgj7}\ba
  |K_8|& =\xl| \int \na P\cdot (1+ \ln (e+|w|^2))wdx\xr|\\
  &\le    \int \n^{\ga-1/2} (1+ \ln (e+|w|^2))\n^{1/2}|\div w|dx+\xl| \int  \n^\ga \na(1+ \ln (e+|w|^2))\cdot wdx\xr|\\
    &\le    C\int \n^{2\ga-1}  \ln^2 (e+|w|^2)dx+\int \n |\div w|^2dx+2\xl| \int  \n^\ga \frac{w\cdot\na w\cdot w}{(e+|w|^2)}dx\xr|\\
    &\le    C\int \n^{2\ga-1}  \ln^2 (e+|w|^2)dx+C\int \n |\na  w|^2dx\\
    &\le    C+C\|\na \n^{\ga/2}\|_{L^2}^2+C\int \n |\na  w|^2dx,\ea\ee
where in the last inequality one has used  the following fact
\be\label{mvgj8}\notag\ba
   \int \n^{2\ga-1}  \ln^2 (e+|w|^2)dx&\le C\int \n^{5\ga/3}  dx+C\int \n   \ln^{\frac{10\ga-6}{3-\ga}} (e+|w|^2)dx\\
   &\le C\|\n\|_{L^\ga}^{2\ga/3}\xl(\|\n\|_{L^1}^\ga+\|\na \n^{\ga/2}\|_{L^2}^2\xr)+C\int \n dx +  C\int \n |w|^2 dx\\
   &\le C+C\|\na \n^{\ga/2}\|_{L^2}^2\ea\ee
owing to Sobolev inequality and \eqref{imszqbb30}.

Finally, the term $K_9$ can be handled  as follows:
\be\label{mvgj9}\ba
   K_9 & =-\ve^{3/2}\int \n |w|^3u\cdot   (1+ \ln (e+|w|^2))wdx \\
  &=-\ve^{3/2}\int \n (1+ \ln (e+|w|^2))|w|^3|u|^2  dx-\mu\ve^{3/2}\int \n (1+ \ln (e+|w|^2))|w|^3u\cdot \na\log\n dx\\
    &\le-\frac{1}{2}\ve^{3/2}\int \n (1+ \ln (e+|w|^2))|w|^3|u|^2  dx+ C\ve^{3/2}\int \n (1+ \ln (e+|w|^2))|w|^3|\na\log\n|^2 dx.\\    \ea\ee
The second term of the right hand of \eqref{mvgj9} holds that for any $0<\beta<\frac{1}{5}$,
\be\label{mvgj10}\ba
    & \ve^{3/2}\int \n (1+ \ln (e+|w|^2))|w|^3|\na\log\n|^2 dx \\
    &\le C \ve^{3/2}\int \n (|w|^3 + |w|^{3+\beta} ) \n^{-1}|\na v|^2 dx  \\
        &\le C \ve^{3/2}\int \n  |w|^5dx + C \ve^{3/2} \int v^{-3}|\na v|^5 dx+C \ve^{3/2} \int \n\n^{-\frac{5}{2-\beta}}|\na v|^{\frac{10}{2-\beta}} dx\\
        &\le C \ve^{3/2}\int \n  |w|^5 dx+ C \ve^{3/2} \int \xl( \n^{-4}+v^{-2}|\na v|^6 \xr)dx+C \ve^{3/2} \int \n^{-\frac{4+3\beta}{6-3\beta}}\xl( v^{-1/3}|\na v|\xr)^{\frac{10}{2-\beta}} dx\\
          &\le C+C \ve^{3/2}\int \n  |w|^5 dx+ C \ve  \int v^{-2}|\na v|^6  dx+C \ve^{3/2} \int \n^{-\frac{4+3\beta}{1-3\beta}} dx\\
          &\le C + C\ve^{3/2}\int \n  |w|^5dx+ C \ve \int v^{-2}|\na v|^6  dx+C \ve^{3/2} \int \n^{-p_0} dx,  \ea\ee
          where one has used  \eqref{imszqbb30} and \eqref{imszqbb5'}.

Submitting \eqref{mvgj5}, \eqref{mvgj4}--\eqref{mvgj10} into \eqref{mvgj1} yields that
\be\label{mvgj11}\ba \notag
 &\frac{1}{2}\frac{d}{dt} \int \n (e+|w|^2)\ln (e+|w|^2)dx +\int \n  \ln (e+|w|^2) \xl(2(\nu-\mu)|\mathcal{D} w|^2+ \sqrt{\ve}\ |\na w|^2\xr)dx\\
 &\le C+C\|\na \n^{\ga/2}\|_{L^2}^2+C\int \n |\na u|^2dx+C\mu^2\int \n |\na^2\log\n|^2dx+ C\ve^{3/2}\int \n  |w|^5dx+ C \ve \int v^{-2}|\na v|^6  dx.\ea\ee
Integrating the upper inequality over $[0,T]$, one obtains \eqref{mvgj0} after suing  \eqref{imszqbb30}, \eqref{imszqbb5'}, and  \eqref{MVinitial}.
 The proof of Lemma \ref{imslem2} is completed. \hfill $\Box$

\emph{Proofs of Theorem  \ref{newth}}: With the   energy estimate, the BD entropy inequality, and the Mellet-Vasseur type estimate obtained in Lemmas \ref{imslem1}--\ref{imslem2} in hand, following the similar compactness arguments as in Section 3 (see also those in \cite{li04,AS2015,AS2017}), one can perform   the  limit progress $\ve\rightarrow 0^+$ to the smooth approximation solutions and thus complete the proof of Theorem  \ref{newth}. We omit the details here.  \hfill $\Box$

\end{document}